\documentclass[]{amsart}

\usepackage{amssymb}
\usepackage{amsmath}
\usepackage{amsthm}
\usepackage{amscd}
\usepackage{graphicx}

\long\def\comment#1\endcomment{}
\usepackage{color}

\newcommand{\vk}{{van Kampen}\ }
\newcommand{\dec}{{criss-cross}\ }

\newcommand{\ertc}{{$\eta$--relator-tied component}}
\newcommand{\ertg}{{$\eta$--relator-tied geodesic}}
\newcommand{\eprt}{{$\eta' $--relator-tied}}

\newcommand{\erts}{{$\eta$--relator-tied}\ }

\newcommand{\ertgs}{{$\eta$--relator-tied geodesic}\ }

\newcommand{\eprtcs}{{$\eta' $--relator-tied component}\ }

\newcommand{\erc}{{$\eta$--relator-covered}}
\newcommand{\ercc}{{$\eta$--relator-covered component}}

\newcommand{\erccs}{{$\eta$--relator-covered component}\ }

\newcommand{\lqq}{\leqslant}
\newcommand{\gq}{\geqslant}
\newcommand {\iv}{^{-1}}

\newcommand{\vf}{\varphi}

\def\trian{\delta}

\def\<{\left\langle}
\def\>{\right\rangle}

\newcommand {\pgot}{{\mathfrak p}}
\newcommand {\rgot}{{\mathfrak r}}

\newcommand {\dgot}{{\mathfrak d}}

\newcommand{\me}{\medskip}

\newcommand {\p}{\mathfrak p} 
\newcommand {\pg}{\mathfrak g} 
\newcommand{\la}{\langle}
\newcommand{\ra}{\rangle}

\newcommand{\fh}{\mathfrak{h}}
\def\fm{{\mathfrak m}}

\newtheorem{thm}{Theorem}[section]
\newtheorem{prop}[thm]{Proposition}
\newtheorem{cor}[thm]{Corollary}
\newtheorem{conj}[thm]{Conjecture}
\newtheorem{qn}[thm]{Question}
\newtheorem{lem}[thm]{Lemma}
\newtheorem{cvn}[thm]{Convention}

\newtheorem{notation}[thm]{Notation}

\theoremstyle{definition}
\newtheorem{defn}[thm]{Definition}
\newtheorem{convention}[thm]{Convention}

\theoremstyle{remark}
\newtheorem{rmk}[thm]{Remark}

\newtheorem{exa}[thm]{Example}

\def\square{\hfill${\vcenter{\vbox{\hrule height.4pt \hbox{\vrule
width.4pt height7pt \kern7pt \vrule width.4pt} \hrule
height.4pt}}}$}

\newcommand{\tsh}[1]{\left\{\kern-.9ex\left\{#1\right\}\kern-.9ex\right\}}

\def\R{{\mathbb R}}

\def\Z{{\mathbb Z}}
\def\N{{\mathbb N}}

\def\C{{\mathbb C}}

\def\dd{{\mathcal D}}

\def\calc{{\mathcal C}}

\def\G{{\mathcal G}}

\def\fg{{\mathfrak g}}

\def\fu{{\mathfrak u}}
\def\fv{{\mathfrak v}}

\def\dist{{\rm{dist}}}

\definecolor{darkgreen}{cmyk}{1,0,1,.2}
\definecolor{Myb}{rgb}{0.9,0,0.9}

\newcommand\intersect\cap
\newcommand\infinity\infty
\newcommand\wt\widetilde

\newcommand\inject\hookrightarrow
\newcommand\union\cup

\newcommand\join\Lambda
\newcommand\cross\times

\newcommand\lub\vee
\newcommand\glb\wedge

\renewcommand\paragraph[1]{\medskip\textbf{#1} }

\begin{document}
\title[Geometry of small cancellation groups, RD and quasi-homomorphisms]{Geometry of infinitely presented small cancellation groups, Rapid Decay and quasi-homomorphisms}
\author{Goulnara Arzhantseva}\thanks{The first author was supported in part by the ERC grant ANALYTIC no. 259527, and by the Swiss NSF, under Sinergia grant CRSI22-130435.}
\address{University of Vienna,
Faculty of Mathematics, Nordbergstra${\ss}$e 15, 1090 Wien, Austria}
\email{goulnara.arzhantseva@univie.ac.at}
\author{Cornelia Dru\c{t}u}\thanks{ The second author was supported in part by
the EPSRC grant ``Geometric and analytic aspects of infinite groups" and by the project ANR Blanc ANR-10-BLAN 0116, acronym GGAA}
\address{Mathematical Institute,
24-29 St Giles,
Oxford OX1 3LB,
United Kingdom.}
\email{drutu@maths.ox.ac.uk}
\date{\today}

\begin{abstract}
We study the geometry of infinitely presented groups satisfying the small cancelation condition $C'(1/8)$, and define a standard decomposition (called the \dec decomposition) for the elements of such groups.  We use it
to prove the Rapid Decay property for groups $G$ with the stronger small cancelation property $C'(1/10)$.
As a consequence, the Metric Approximation Property holds for  the reduced $C^*$--algebra $C^*_r(G)$ and for the Fourier algebra $A(G)$ of the group $G$. Our method further implies that the kernel of the comparison map between the bounded and the usual group cohomology in degree $2$ has a basis of power continuum.

The present work can be viewed as a first non-trivial step towards a systematic investigation of direct limits of hyperbolic groups.

\end{abstract}
\maketitle

\me
\section{Introduction}

The construction of infinite finitely generated groups with unusual properties, the so-called ``infinite monsters'', plays a major part in geometric group theory. The purpose of such groups is to test the robustness and the level of generality of various long-standing conjectures. A main source of infinite monsters is the class of infinitely presented direct limits of Gromov hyperbolic groups, and within it the class of small cancellation groups, in their numerous variants: classical, of Olshanskii type, graphical, geometric etc. This is the case for the Tarski monster groups constructed by Olshanskii \cite{Olsh:book}, for groups with non-homeomorphic asymptotic cones such as the Thomas-Velickovic examples \cite{ThomasVelickovic} and the Dru\c{t}u-Sapir  examples \cite{DrutuSapir:TreeGraded}, and for the Gromov monster groups containing expander families of graphs in their Cayley graphs \cite{Gromov:random,ArzhDelzant}.

As mentioned, infinite monsters are usually designed to be counter-examples to various statements, and it is rather challenging to obtain positive results about them.  The few known such results are algebraic or geometric in nature. Analytic properties have not yet been investigated systematically or with any attempt at generality so far. Their study did not go beyond that of a few historically foundational examples, such as the free Burnside groups of sufficiently large odd exponent, or the Tarski monster groups, relevant to the property of (non)-amenability, or the Gromov monster groups, relevant among others to Kazhdan's property (T), etc.
At the same time, the (counter)-examples required by K-theory, operator algebra and topology and the existing constructions of infinite monster groups do incite to such a research.
The present work can be viewed as a first step in this direction.

In this paper,  we prove two general results on infinitely presented small cancellation groups.  The first result yields several analytic properties of such groups, in particular the  property of Rapid Decay, see Section \ref{subsec:first} and Theorem \ref{thm:main}. The second result provides a new way of constructing an abundance of quasi-homomorphisms in these groups, see Section \ref{subsec:second} and Theorem \ref{thm:qm}. Our methods very likely generalize to larger classes.

Note that there are known examples of infinitely presented small cancellation groups that appear as subgroups of Gromov hyperbolic groups \cite[Theorem 3.1 and 4.3]{KapovichWise}. Therefore these groups have the property RD, the Haagerup property, and all their consequences; they also have a large space of quasi-homomorphisms by work of Bestvina-Fujiwara \cite{BestvinaFujiwara}. However, these are quite specific examples.

\subsection{Main technical tool: \dec decompositions}
For several analytic and geometric group properties, including the two discussed in this paper, it is crucial to understand if the group elements possess  ``standard decompositions'' into products of  certain ``elementary'' parts. These decompositions are required to behave well with respect to the group structure, that is, when one considers triples of elements $g_1,g_2$ and $g_3$ such that $g_1g_2=g_3\, $.

The main technical result of our paper is the construction of such a decomposition for elements of finitely generated groups defined by infinite presentations with the small cancellation condition $C'(\lambda)$, for $\lambda \lqq \frac{1}{8}$ (so-called $C'(1/8)$--groups). More precisely,
given a pair of vertices in a Cayley graph of such a group, we obtain a detailed description of a set containing all the geodesics between the two vertices, see Theorem \ref{lem:intermediate}.
These sets, in many ways, play the role of the convex hulls used in \cite{RamaggeRobSteg} to show the property of Rapid Decay for groups acting on buildings. The existence of such sets allows us to introduce a uniquely defined decomposition, called the \emph{\dec decomposition}, of the elements of the given $C'(1/8)$--group,  see Section \ref{sec:stdec}.

It is worth noticing that our approach differs and cannot be deduced from the Rips-Sela canonical representatives \cite{Sela:uniform,RipsSela:canonical} in finitely presented small cancellation groups. Indeed, the Rips-Sela construction gives an equivariant choice of quasi-geodesic paths between pairs of vertices (with a view to reduce solving equations in a small cancellation group, or more generally in a hyperbolic group, to solving equations in a free group). Their arguments are based on the existence of central points for geodesic triangles, granted by finite presentation only. The similarity in method between the Rips-Sela approach and ours does not go beyond the common use of the geometry of geodesic bigons \cite[Theorem 5.1]{RipsSela:canonical}.

\subsection{First result: property of Rapid Decay}\label{subsec:first} The property of Rapid Decay (or property RD; see Section \ref{section:rd} for the definition) can be seen as a non-commutative generalization of the density of the set of smooth functions on a torus $\mathbb{T}^n$ inside the algebra of continuous functions on $\mathbb{T}^n\, ;$  it requires the existence of an analogue of the Schwartz space inside the reduced $C^*$--algebra of a group, see Example \ref{exa:zn}.

Property RD was proved for free groups by Haagerup in \cite{Haagerup:RD_F}. It was formally defined and established for several classes of groups by Jolissaint in \cite{Jolissaint:defRD}. It easily follows from the definition of property RD that polynomial growth implies this property \cite{Jolissaint:defRD}. Other examples of groups known to possess RD include (the list below is not exhaustive):

\begin{itemize}

\item Gromov hyperbolic groups \cite{Harpe:RDforhyp,Jolissaint:defRD};

\item groups hyperbolic relatively to subgroups with property RD \cite{DrutuSapir:RD};

\item groups acting on certain buildings or symmetric spaces of rank two \cite{RamaggeRobSteg,Lafforgue:RD};
groups acting on products of buildings or symmetric spaces of rank one or two \cite{Chatterji:RDproducts};

\item mapping class groups of surfaces \cite{BehrstockMinsky:RD};

\item large type Artin groups \cite{CiobanuHoltRees:RD}.

\end{itemize}

The only known obstruction to property RD is due to Jolissaint \cite{Jolissaint:defRD}: a group that contains an amenable subgroup
of super-polynomial growth (with respect to the word length metric of the ambient group) does not have property RD. For example, every non-uniform irreducible lattice in a semi-simple group of rank at least two fails to have property RD.

Despite the fact that property RD has many applications (see Section \ref{sect:appl} for a sample), the full extent of the class of groups satisfying RD is yet to discover. Moreover, with a few exceptions, the examples of finitely generated groups with property RD known up to now are finitely presented.
On the other hand, most relevant counterexamples in
$K$-theory and topology, such as infinite simple groups with Kazhdan's property (T) or Gromov's monster groups, are infinitely presented.

In this paper, as a \emph{first application} of the \dec decomposition, we prove that infinitely presented $C'(1/10)$--groups have the
property of Rapid Decay.

\begin{thm}\label{thm:main}
Every finitely generated group defined by an infinite presentation satisfying the small cancelation condition $C'(1/10)$ has the property of Rapid Decay.
\end{thm}

Since finitely presented $C'(1/10)$--groups
are Gromov hyperbolic, a group as in Theorem \ref{thm:main} is a direct limit of Gromov hyperbolic groups, hence of groups with property RD. However,  property RD is not preserved under taking direct limits: there exist elementary amenable groups of exponential growth which are direct limits of
Gromov hyperbolic groups \cite{OOS}.

\begin{qn}\label{qn:RDlimit}
Does the class of direct limits of hyperbolic groups satisfying the RD property contain all $C'(1/6)$--groups? Does it contain the free Burnside groups of sufficiently large odd exponent, Tarski monsters, Gromov's monsters?
\end{qn}

Theorem \ref{thm:main} is meaningful in the context of the Baum-Connes conjectures. Namely,  the existing proofs of various Baum-Connes conjectures
encounter two types of obstructions:
\begin{itemize}
  \item Kazhdan's property (T) and its strengthened versions \cite{Lafforgue:Trenforce}. These are obstructions to a certain strategy of proof.
  \item An expander family of graphs coarsely embedded in the Cayley graph of a group.  Such a group, known as a Gromov monster \cite{Gromov:random, ArzhDelzant},
   does not satisfy the Baum-Connes conjecture with coefficients \cite{HigsonLafforgueSk}, see also \cite{WYu:index1,WYu:index1}. The original Baum-Connes conjecture (without coefficients or, equivalently, with complex coefficients) is still open for Gromov monsters.
\end{itemize}

In many cases, the property (T) obstruction can be overcome: V. Lafforgue proved  that  the property of Rapid Decay combined with a ``good'' action on a CAT(0)--space implies the Baum-Connes conjecture without coefficients \cite{Lafforgue:Baum-Connes}. This applies to many groups with property (T) such as certain hyperbolic groups and uniform lattices in $SL(3, \R )$.

Our methods may contribute to overcome the second type of obstruction in the same way as the first.
Indeed, Gromov's monster groups have graphical small cancellation presentations and our results possibly extend to prove property RD for such groups.
Note that a recent result of Osajda and the first author \cite{ArzhOsajda:Haagerup} implies that all infinitely presented $C'(1/6)$--groups have the Haagerup property (see also Section \ref{sec:applapprox}), hence they satisfy the Baum-Connes conjecture with coefficients \cite{HigsonKasparov:HaagerupBC}. However, small cancellation groups in general cannot satisfy the Haagerup property: Gromov's monsters do not even admit a coarse embedding in a Hilbert space.

\subsection{Applications of the Rapid Decay property}\label{sect:appl}
Combined with appropriate analytic and geometric properties, property RD yields non-trivial applications to non-commu\-ta\-tive geometry, harmonic analysis, growth properties and quantum dynamics. We provide more details in Section \ref{sec:mainappl}, here we overview them briefly.

The first application concerns finite-dimensional approximations of $C^*$--algebras and other operator algebras of discrete groups. The study of such approximations was initiated by Grothendick \cite{Grothendieck:book} and received a further impetus through the Banach space counter-examples of Enflo \cite{ Enflo:counterexAP} (see also \cite[Section 2.d]{LindenstraussTzaf1} and \cite[SectionI.g]{LindenstraussTzaf2}),  and the $C^*$--algebra counter-examples of Szankowski  \cite{Szankowski:counterex}, who emphasized that such properties are far from generic. We refer to \cite{CCJJV,BrownOzawa:book,LafforguedelaSalle} and references therein for details.

 For reduced $C^*$--algebras of infinite groups, approximation properties relate to various kinds of amenability (e.g. the classical amenability of von Neumann, weak amenability in the sense of Cowling-Haagerup, $C^*$-exactness due to Kirchberg-Wassermann, etc). For instance, the reduced $C^*$--algebra of a group is nuclear if and only if that group is amenable \cite{Lance:AmenabNucl}.

The non-abelian free groups are the first examples of groups with a reduced $C^*$--algebra which, while non-nuclear, has the Metric Approximation Property \cite{Haagerup:RD_F}. This latter property means that the identity map on the $C^*$--algebra can be approximated in the strong topology (also called the point-norm topology) by a net of finite rank contractions \cite[Definition 1.e.11]{LindenstraussTzaf1}. The arguments in \cite{Haagerup:RD_F} inspired the result that we  will use: if $G$ has the RD property with respect to a conditionally negative definite length function then the Metric Approximation Property holds for the reduced $C^*$--algebra $C^*_r(G)$ and for the Fourier algebra $A(G)$ of $G$ \cite{Haagerup:RD_F,JolissaintValette:RDapp,BrodNiblo:approximation}.

In our context, it is therefore necessary to find a length function on a $C'(1/10)$--group that is conditionally negative definite and bi-Lipschitz equivalent to a word length function. We explain in Section \ref{sec:mainappl} how this can be deduced from the work of Wise \cite{Wise:cancell} and the \dec decomposition described in Theorem \ref{lem:intermediate}. We refer the reader to~\cite{ArzhOsajda:Haagerup} for a stronger result and a different approach.

Theorem \ref{thm:main} then implies

\begin{cor}\label{cor:map}
Let $G$ be a finitely generated group given by an infinite presentation satisfying the small cancelation condition $C'(1/10)$.
Then the reduced $C^*$--algebra $C^*_r(G)$ and the Fourier algebra $A(G)$ of $G$ have the Metric Approximation Property.
\end{cor}

Our \dec decomposition can be used to get other strong approximation results for the algebras $C^*_r(G)$ and $A(G)$ of $C'(1/10)$--groups, in particular in connection to the question that we formulate below. We postpone these considerations to an upcoming paper.

\begin{conj}\label{conj:wa}
Let $G$ be a finitely generated group defined by an infinite presentation satisfying the small cancelation condition $C'(1/10)$.
Then $G$ is weakly amenable and $C^*$-exact (equivalently, it has Guoliang Yu's property A).
\end{conj}

The methods developed in the proof of Theorem \ref{thm:main} suggest an affirmative answer to Conjecture \ref{conj:wa}. Namely, an intermediate step in the proof of the Metric Approximation Property of the Fourier algebra is that $A(G)$ has an approximate identity bounded in the multiplier norm \cite{JolissaintValette:RDapp, BrodNiblo:approximation}, cf. \cite{Haagerup:RD_F,CanniereHaagerup:Fourier}.
If replaced by a completely bounded multiplier norm, the result would mean that $G$ is weakly amenable. Equivalently, the reduced $C^*$--algebra of $G$ would satisfy the completely bounded approximation property \cite{HaagerupK:approximation}.  This approximation property implies the $C^*$-exactness. For general direct limits of hyperbolic groups, the conjecture does not hold  as  Gromov's monster group containing an expander in its Cayley graph is not $C^*$-exact, and hence, it is not weakly amenable.
Finitely presented $C'(1/6)$ small cancellation groups are both weakly amenable and  $C^*$-exact as these are known for Gromov hyperbolic groups \cite{Ozawa:wa,Adams:exact}.

 The next application concerns the growth series of small cancellation groups. It follows by combining our property RD result with
 that of \cite{GrigorchukNagnibeda}, see Section~\ref{sec:mainappl}.

\begin{cor}\label{cor:opg}
Let $G$ be an infinitely presented finitely generated $C'(1/10)$--group.
Then the radius of convergence of the operator growth series equals the square root of the radius of convergence of the standard growth series.
\end{cor}

Finally, property RD is relevant to ergodic theorems in the setting of quantum dynamics.
We give a brief overview on this in Section \ref{sec:quantum}.

\subsection{Second result: the bounded versus the usual cohomology}\label{subsec:second}

A \emph{second application} of  the \dec decomposition is  that infinitely presented $C' \left( 1/12 \right)$--groups are
rich in quasi-homomorphisms (see Section \ref{sec:qm} for definitions). This has an immediate impact on the bounded cohomology of such groups \cite{Gromov:Volume}.

\begin{thm}\label{thm:bch}
Let $G$ be a  finitely generated group defined by an infinite presentation satisfying the small cancelation condition $C'(1/12)$.
Then the kernel of the comparison map between the second bounded and the usual group cohomology $$
H^2_b(G) \to H^2 (G)\, ,
$$
is an infinite dimensional real vector space, with a basis of power continuum.
\end{thm}

The above kernel can be identified with the real vector space $\widetilde{QH}(G)$ of all quasi-homomorphisms modulo near-homomorphisms
(where by a near-homomorphism we mean a function $\fh :G \to \R$ that differs from a homomorphism by a bounded function).
The computation of $\widetilde{QH}(G)$ is therefore important and it has been done up to now for various classes of groups.

Groups that have a certain type of action on a hyperbolic space (in particular, subgroups of relatively hyperbolic groups, mapping class groups, etc.) have $\widetilde{QH}(G)$ infinite dimensional,
with a basis of power continuum. This was proved by Brooks for non-abelian free groups \cite{Brooks:bc} and by Brooks and Series for non-amenable surface groups. In \cite{Gromov(1987)} Gromov stated that all non-elementary hyperbolic groups have non-trivial second bounded cohomology. Epstein and Fujiwara proved  that in fact for all non-elementary
hyperbolic groups $\widetilde{QH}(G)$ has a basis of power continuum \cite{Epstein-Fujiwara}. Later this result was extended to other types of groups acting on hyperbolic spaces
and to their non-elementary subgroups \cite{Fujiwara:BC1, Fujiwara:BC2}; in particular, to subgroups of mapping class groups of surfaces \cite{BestvinaFujiwara}. See also the survey of Fujiwara \cite{Fujiwara:hb} and references therein. The same result was further extended to groups with free hyperbolically embedded subgroups by Hull and Osin \cite{HullOsin:qcocycles}.

At the other end, large arithmetic lattices have trivial $\widetilde{QH}(G)$ \cite{Monod:book}.

Our approach differs from all the previous ones showing that $\widetilde{QH}(G)$ has a large basis. Indeed, all the previously known arguments rely essentially on the existence of free subgroups. We do not require the existence of such subgroups, and a potential extension of our methods may apply to groups satisfying other small cancellation conditions such as the Olshanskii small cancellation, in particular, to the free Burnside groups of sufficiently large odd exponent or to the Tarski monsters.

The following result is another immediate consequence of our theorem above.

\begin{cor}\label{cor:notb}
Let $G$ be a finitely generated group given by an infinite presentation satisfying the small cancelation condition $C'(1/12)$.
Then $G$ is not boundedly generated\footnote{A group is \emph{boundedly generated} if it can be expressed (as a set) as a finite product of cyclic subgroups.}.
\end{cor}

\subsection{Plan of the paper} The paper is organized as follows. Sections \ref{section:rd} and \ref{section:cancell}  give
preliminary information on the property of Rapid Decay and small cancelation groups. In Section \ref{sec:stdec}, we  describe the \dec decomposition of elements in infinitely presented small cancelation groups. We believe this description is of independent interest and
it can be applied to get further results on such groups.
In Section \ref{sec:main}, we use this \dec decomposition to prove our Theorem \ref{thm:main}. In Section \ref{sec:mainappl} we explain several applications of  Theorem \ref{thm:main}.
In Section \ref{sec:qm} we focus on quasi-homomorphisms of $C'(1/12)$-small cancellation groups and prove Theorem~\ref{thm:bch}.

\medskip

\emph{Acknowledgments.} The work on this paper was carried out during visits of the first author to the University of Oxford and of the second author to the University of Vienna. We thank these institutions for their support and hospitality. The first author is also grateful to the ETH Z\"urich and to the CRM in Barcelona for their hospitality during the completion of the paper.

We thank Thomas Delzant and Denis Osin for useful comments and corrections.

\section{Word metrics and the Rapid Decay Property}\label{section:rd}

Let $G$ be a group generated by a finite subset $A=\{ a_1^{\pm 1},\dots,a_m^{\pm 1}\}$ not containing the neutral element $1$.

Given a word $w$ in the alphabet $A$, we let $|w|_A$ denote its length. For an element $g\in G$ we let $|g|_A$ denote the minimal length of a (reduced) word $w$ in $A$ representing $g$. The function $G \to \N$ defined by $g\mapsto |g|_A$ is called a \emph{word length function}, in particular it is a length function in the usual sense (see for instance Definition 1.1.1 in \cite{Jolissaint:defRD}). It allows to define a left-invariant metric on $G$ by  $\dist_A (x,y)= |x^{-1}y|_A$ for $x,y\in G$, called a \emph{word metric on}~$G$:

For all the problems treated in this paper the choice of the finite generating set $A$ is irrelevant, thus from now on we tacitly assume that such a set is fixed for every group that we consider, and defines length function and metric as above, and we omit the subscripts from the notation.

We denote by
$
B (x, r)=\{
y\in G \mid \dist (x,y)\leqslant r\}
$, the closed ball of radius $r$ centered at $x$, and by
$S (x, r)=\{ y\in G \mid \dist (x,y)= r\}
$ the sphere of radius
$r$ centered at $x$.

We begin by recalling the analytic version of property RD, in order to emphasize its interest, then we formulate other equivalent definitions that are usually easier to prove. Throughout the paper we denote by $\ell^2 (G)$ the Hilbert space of square-summable $\mathbb C$-valued functions on $G$, and by $\| f \|$ the $\ell^2$-norm of a function $f$.

The {\it group algebra of $G$}, denoted by $\C G$, is the set of complex valued functions with finite support on $G$, that is the set of formal
linear combinations of elements of $G$ with complex coefficients. We
denote by $\R_+ G$ its subset consisting of functions taking values
in $\R_+$. The action of $G$ by  left-translation on the space $\ell^2(G)$ extends by linearity, using the convolution, to a faithful action of the group algebra $\C G$ on $\ell^2(G)$:   $f\ast g(z) = \sum_{x\in G} f(x)g(x^{-1}z)$.

We can thus identify $\C G$ with a linear subspace in the space of bounded operators $B\left( \ell^2(G) \right)$. Its closure in the operator norm is denoted by $C^*_r (G)$ and it is called the \emph{reduced $C^*$--algebra of }$G$.

Note that $C^*_r (G)$ is embedded in $\ell^2(G)$ by the map $T\mapsto T(\delta_1)$, where $\delta_1$ is the characteristic function of the singleton set $\{ 1\}\, ,$ with $1$ the identity element in $G$.

\begin{exa}\label{exa:zn}
The particular case of $G=\Z^n$ helps to  understand the general idea of the Rapid Decay property.

The reduced $C^*$--algebra $C^* (\Z^n )$ is isomorphic to the $C^*$--algebra of continuous functions on the $n$--dimensional torus $C(\mathbb{T}^n)$. The latter can be identified to the former \emph{via} the Fourier transform.

On the other hand, $C(\mathbb{T}^n)$ contains a sub-algebra that is dense and full, the algebra of the smooth functions $C^\infty (\mathbb{T}^n)\, $. It remains to note that the Fourier transform of a function in $C^\infty (\mathbb{T}^n)\, $ is a function that ``decays rapidly'', since its product with any power function is still square-summable.
\end{exa}

A generalization of the property above for an arbitrary finitely generated group is as follows.

For every $s\in \R$, {\it the
Sobolev space of order $s$ on the group $G$} is the set $H^s
(G)$ of functions $\phi$ on $G$ such that the function $g\mapsto (1+|g|)^s\phi (g)$
is in $\ell^2 (G)$.

The space of {\it{rapidly decreasing functions on
$G$}} is the set $H^\infty (G)=\bigcap_{s\in
\R }H^s (G)$.

\begin{defn}[Rapid Decay property]\label{def:rd0}
The group $G$ has the \emph{Ra\-pid Decay property} (or \emph{property RD}) if $H^\infty (G)$ is contained in $C_r^*(G)$, seen as a subspace of $\ell^2 (G)\, $.

Equivalently, if the inclusion of $\C G$ into $C_r^*(G)$ extends to a continuous inclusion of $H^\infty(G)$ into $C_r^*(G)$.
\end{defn}

The general definition of the RD property requires that for some length function on $G$ the corresponding space of rapidly decreasing functions is contained in the reduced $C^*$--algebra. Still, for finitely generated groups, this is equivalent to the property for one (for every) word length function (see for instance \cite{DrutuSapir:RD} for details).

The significance of property RD is emphasized by a consequence of it which goes back to Swan and Karoubi stating that if $H^\infty (G)$ is contained in $C_r^*(G)$ then the inclusion induces an isomorphism of
$K$--groups $K_i \left( H^\infty (G) \right)$ with $K_i \left( C_r^*(G) \right)$, for $i = 0, 1$.
 This is a main ingredient in the proof due to Connes and Moscovici of the Novikov conjecture for
Gromov hyperbolic groups  \cite{ConnesMoscovici}.

An equivalent definition of the Rapid Decay property is as follows.

\begin{defn}[Rapid Decay property]\label{def:rd}

A finitely generated group $G$ has the \emph{property RD} if there exists a polynomial $P$ such that for every $R>0$, every function $f\in \C G$ vanishing outside the ball $B(1,R)$, and every $g\in l^2(G)$, the following inequality holds
\begin{equation}\label{ineq0}
\| f*g \| \leqslant P(R)\, \| f \| \cdot \| g \|\, .
\end{equation}
\end{defn}

In fact, (\ref{ineq0}) holds whenever it is satisfied by functions with positive values and finite support, as shown by the next lemma.

For a function $f\in l^2(G)$ and a constant $p \geqslant 0$, $f_p$
denotes the function which coincides with $f$ on $S(1,p)$ and which
vanishes outside $S(1,p)$.

\begin{lem}\label{equiv1} Let $G$ be a finitely generated group. The following statements are equivalent:
\begin{itemize}
    \item[(i)] The group $G$ has property RD.
\item[(ii)] There exists a polynomial $P$
such that for every
    $r, R \gq 0$, every $p\in [| r-R | \, ,\, r+R]$
    every $f\in \R_+G$
    with support in
    $S(1,R)$, and every $g\in \R_+ G$
    with support in $S(1, r)$,
\begin{equation}\label{ineq}
    \| (f*g)_p \| \leqslant P(R)\, \| f \| \cdot \|g \|\, .
\end{equation}
\end{itemize}
\end{lem}

The details and proofs of the above equivalences, see the proof of Theorem  5 \cite[$\S$ III.5.$\alpha $]{Connes:noncommutative},  \cite{ChatterjiRuane}, and Lemma 2.7 in \cite{DrutuSapir:RD}.

For further information on property RD and its applications we refer to \cite{Haagerup:RD_F, Jolissaint:defRD, Connes:noncommutative, Lafforgue:RD,ChatterjiRuane,DrutuSapir:RD}.

\section{Preliminaries on infinite small cancelation presentations}\label{section:cancell}

A set of words $R$ in the alphabet $A$ is said to be \emph{symmetrized} if it contains $r\iv$ and all the cyclic permutations of $r$ and $r\iv$, whenever $r\in R$.
Without loss of generality we always assume that the set of group relators is symmetrized and that all relators $r\in R$ are reduced words in the alphabet $A$.

We focus on groups with \emph{infinite presentations},
\begin{equation}\label{eq:pres}
G=\la A \mid r_1,\ldots, r_k,\ldots  \ra\, ,
\end{equation}
defined by a symmetrized family $R$ of relators consisting of an infinite sequence of relators $r_1,\ldots , r_k,\ldots $.

We denote by $R_k$ the set $\{r_1, \ldots,  r_k\, \} $ and by $G_k$ the finitely presented group
\begin{equation}\label{eq:presK}
G_k= \la A \mid R_k \rangle\,= \la A \mid r_1, \ldots,  r_k \ra\, .
\end{equation}

For two words $u,v$ we write $u\sqsubset v$ when $u$ is a subword of $v$.
Let $\eta$ be a constant in $\left(0, \frac{1}{2}\right]$. If in the preceding we have moreover that
$$\eta |v| \lqq |u| \lqq \frac{1}{2} |v|$$ then we use the notation $u\sqsubset_\eta v$.
We write $u\sqsubset R$ if there exists $v\in R$ such that $u\sqsubset v$. Similarly, for $\sqsubset$ replaced by  $\sqsubset_\eta.$

\begin{notation}\label{notat:seta}
We denote by $S(R)$ the set of words $u$ such that  $u\sqsubset R$ and by $S^\eta (R)$ the set of words $u$ such that  $u\sqsubset_\eta R.$
\end{notation}

\begin{defn}[$C'(\lambda)$--condition]\label{def:clambda} Let $\lambda\in\left (0, 1\right)$.
A symmetrized set $R$ of words in the alphabet $A$ is said to satisfy the
{\it $C'(\lambda )$--condition} if the following holds:
\begin{itemize}
\item[(1)] If $u$ is a subword in a word $r\in R$ so
that $|u|\gq \lambda |r|$ then $u$ occurs only once in $r$;
\item[(2)] If $u$ is a subword in two distinct words $r_1,r_2\in
R$ then $|u|< \lambda \min \{|r_1|,|r_2|\}$.
\end{itemize}
We say that a group presentation $\la A \mid R\ra $ {\it satisfies $C'(\lambda )$--condition} if $R$ satisfies that condition.
\end{defn}

Our technical arguments use the language of van Kampen diagrams over a group presentation $\la A \mid R\ra ,$
for more details see \cite{LyndonSchupp} (observe that the classical results below still hold for infinite group presentations).

The boundary of any \vk diagram (cell)  $\Delta$  is
denoted by $\partial\Delta.$

\begin{lem}[Greendlinger \cite{LyndonSchupp}{Ch.V, Thm. 4.4}]\label{lem:Green}
Every reduced \vk diagram $\Delta $ over the presentation \eqref{eq:pres} with small cancelation condition $C' (\lambda )$ for $\lambda \lqq \frac 16$ contains a cell $\Pi$ with $\partial \Pi$ labeled by a relator $r\in R$ such that $\partial \Delta \cap \partial \Pi$ has a connected component of length  $>(1-3\lambda)\,  |r|$.
\end{lem}

\begin{defn}[$n$-gone]\label{def:simplegone}
We call $n$\emph{-gone} in a geodesic metric space a loop obtained by successive concatenation of $n$ geodesics.

We say that the $n$-gone is \emph{simple} if the loop thus obtained is simple, that is, it does not have self-intersections.
\end{defn}

\begin{figure}[htb]
\centering
\includegraphics[scale=0.6]{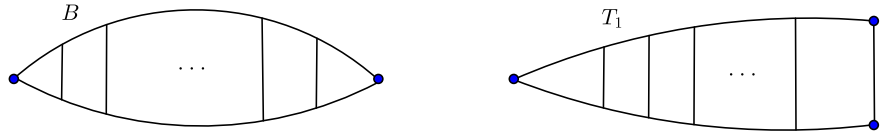}
\caption{Simple bigon $B$ and simple triangle $T_1$.}
\label{fig:bigon}
\end{figure}

\begin{figure}[htb]
\centering
\includegraphics[scale=0.9]{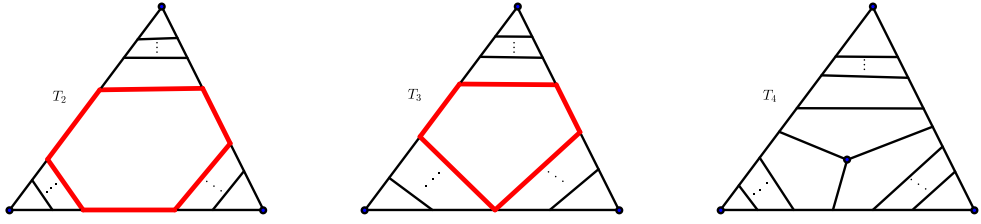}
\caption{Simple triangles  $T_2, T_3,$ and $T_4$.}
\label{fig:triangle}
\end{figure}

\begin{thm}[cf. \cite{Ghys-Harpe(1990)}]\label{thm:classif}
Let $\Delta $ be a reduced \vk diagram over a group presentation $G= \la A \mid R\ra$ satisfying the $C'(\lambda )$--condition, with $\lambda \leqslant \frac 18$.
\begin{enumerate}
  \item\label{bigon} Assume that $\partial \Delta $ is a simple bigon in the Cayley graph of $G$. Then it has the form of the bigon $B$ in Figure \ref{fig:bigon}.
  \item\label{triangle}  Assume that $\partial \Delta $ is a simple triangle in the Cayley graph of $G$.  Then it has one of the forms $T_1, \dots, T_4$ in Figure \ref{fig:bigon} and Figure \ref{fig:triangle}.
\end{enumerate}
\end{thm}

\section{Standard decomposition of elements in small cancellation groups}\label{sec:stdec}

This section is devoted to a thorough analysis of geodesics in Cayley graphs of infinitely presented small cancelation groups, and to the description of a set which, from many points of view, plays the part of the
 convex hull of a two points-set in irreducible buildings.
We show here the main technical result of the paper, Theorem \ref{lem:intermediate}, and its algebraic counterpart Theorem \ref{cor:sdecomp}.

\begin{convention}\label{conv:g18}
Throughout this section $G$ denotes a finitely generated group with a (possibly infinite) presentation $\la A \mid R \ra$ satisfying the $C'(\lambda )$--condition with $\lambda \leqslant \frac{1}{8}$.

We only consider the Cayley graph of $G$ with respect to the fixed (arbitrary) finite generating set $A$, and we omit mentioning $A$ from now on. By ``vertex'' we shall always mean a vertex in that Cayley graph.
\end{convention}

We call \emph{contour} a loop in the Cayley graph of $G$ labeled by a relator $r\in R$. By abuse of notation, given a contour $t$ we denote by $|t|$ its length. Observe that a contour is always a simple loop (a non-trivial self-intersection leads to a contradiction with the small cancelation assumption by the Greendlinger lemma).

By an \emph{arc} we mean a topological arc, that is the image of a topological embedding of an interval into a topological (in particular metric) space.

 For every path $\pgot$ in a metric space, we denote
the initial point of $\pgot$ by $\pgot_-$ and the terminal point of $\pgot $ by $\pgot_+$.
Given two points $x,y$ on a geodesic $\pg$, we denote by $[x,y]$ the sub-geodesic of $\pg$ with endpoints $x,y$.

\begin{lem}\label{lem:geod}
Let $t$ be a contour labeled by a relator $r$ and let $a,b$ be two points on $t$.
\begin{enumerate}
  \item\label{geod1} If one of the two arcs with endpoints $a,b$ have length $<\frac{|t|}{2}$ then that arc is the unique geodesic with endpoints $a,b$ in the Cayley graph.

  \me

  \item\label{geod2} If both arcs with endpoints $a,b$ have length $\frac{|t|}{2}$ then these arcs are the only two geodesics with endpoints $a,b$ in the Cayley graph.

  \item\label{geod3} The intersection of a geodesic with a contour is always composed of only one arc.
\end{enumerate}
\end{lem}

\proof \eqref{geod1} Assume there exists a geodesic joining $a,b$ distinct from that arc. Then they compose at least one non-trivial simple bigon. Consider the minimal \vk diagram $\Delta$ with boundary labeled same as this bigon. Let $u$ be the label of the sub-arc of $t$ and $v$ the label of the sub-arc of the geodesic. According to Lemma \ref{lem:Green}, there exists a cell $\Pi$ labeled by a relator intersecting the boundary $\partial \Delta$ in an arc of length $>1-3\lambda$ of the length of $\partial \Pi$.

Assume first that $\partial \Pi$ does not coincide with $t$. By the small cancelation condition, the arc can have at most $\lambda $ of the length of $\partial \Pi$ in common with the arc labeled by a subword of $r$, hence it has $>1-4\lambda$ of the length of $\partial \Pi$ in common with the arc labeled same as the geodesic. As $\lambda \leqslant \frac{1}{8},$ this contradicts the fact that this is the label of a geodesic.

Now if $\partial \Pi$ coincides with $t,$ then $\partial \Pi$ has at least $\frac{1}{2}-3\lambda$ of its length in common with the arc labeled by $v$. In particular, it follows that $|u| \geqslant |v| > \left(\frac{1}{2}-3\lambda \right) |r|$. Then $\partial \Delta \bigtriangleup \partial \Pi$ composes a new simple bigon with both sides of length at most $3\lambda |r|.$ We apply the argument above to this new bigon, the boundary of the cell provided by Lemma \ref{lem:Green} cannot coincide with $t$ this time, and we obtain a contradiction.

\me

\eqref{geod2} The argument to show that there exists no geodesic joining $a,b$ and which is not entirely contained in $t$ is as above.

\me

\eqref{geod3} It suffices to prove that this intersection is path connected. Indeed, let $\pg$ be a geodesic and let $a,b$ be two points on $\pg \cap t$. The above arguments show that the part of $\pg$ between $a$ and $b$ must be contained in $t$.\endproof

\begin{defn}[Relator-tied geodesics and components]
Let $\pg $ be a geodesic in the Cayley graph of $G$ and let $\eta$ be a number in $(0,1)$.
\begin{enumerate}
  \item $\pg$ is called $\eta$--\emph{relator-tied} if it is covered by sub-geodesics labeled by words in $S^\eta (R)$.
  \item an $\eta$--\emph{relator-tied component} of $\pg$ is a maximal sub-geodesic of $\pg$ that is $\eta$--relator-tied.
\end{enumerate}
\end{defn}

\begin{lem}\label{lem:etacomp}
\begin{enumerate}
  \item\label{comp1} The $\eta$--relator-tied components of a geodesic $\pg $ are disjoint.

  \item\label{comp2} Assume that $\eta \lqq \frac{1}{2}-2\lambda$. If two points $a,b$ are the endpoints of a geodesic $\pg$ with no $\eta$--relator-tied component then $\pg$ is the unique geodesic with endpoints $a,b$.
\end{enumerate}
\end{lem}

\proof Assertion \eqref{comp1} follows by definition, since two distinct $\eta$--relator-tied sub-geodesics that intersect compose a longer $\eta$--relator-tied sub-geodesic.

\eqref{comp2} \quad Any other geodesic $\pg'$ with endpoints $a,b$ and distinct from $\pg$ would compose with $\pg$ simple geodesic bigons, therefore by Theorem \ref{thm:classif}, \eqref{bigon}, $\pg$ would contain a $\left( \frac{1}{2}-2\lambda \right)$--relator-tied component.\endproof

\begin{defn}[$\eta$--compulsory geodesic]
Given $0< \eta \lqq \frac{1}{2}-2\lambda$, a geodesic as in Lemma \ref{lem:etacomp}, \eqref{comp2}, is called an $\eta$--\emph{compulsory geodesic}.
A pair of endpoints $a,b$ of an $\eta$--compulsory geodesic is called an $\eta$--\emph{compulsory pair}.
\end{defn}

We now proceed to analyze the $\eta$--relator-tied components of geodesics.

\begin{figure}[htb]
\centering
\includegraphics[scale=0.8]{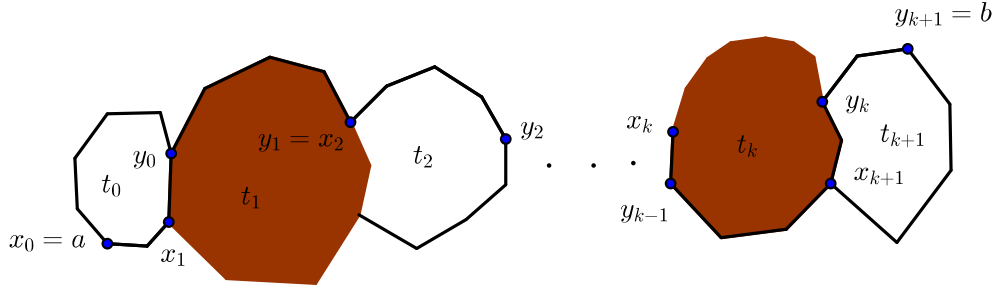}
\caption{An $\eta$--relator-tied geodesic inside a succession of contours.}
\label{fig:RTgeodesics}
\end{figure}

\begin{lem}\label{lem:rt1} Let $\eta \gq 2\lambda$.
 Let $\pg$ be a $\eta$--relator-tied geodesic in the Cayley graph of $G.$ Then there exists a unique sequence of successive vertices $$x_0=a, x_1, y_0, x_2, y_1,\ldots, x_{k+1}, y_k , y_{k+1}=b$$
  such that the sub-geodesics with endpoints $x_i, y_i$ with $i\in \{ 0,1,\ldots,k+1\}$ are labeled by words in $S^\eta (R)$, and are maximal with this property with respect to inclusion (see Figure \ref{fig:RTgeodesics}).
\end{lem}

\proof By hypothesis $\pg \subseteq \bigcup_{i\in S_0} \pg_i$, where $\pg_i$ denotes a sub-geodesic of $\pg$ labeled by a word $u_i\in S^\eta (R)$ and the index set $S_0$ is finite (by compactness).

Without loss of generality, we assume that all the sub-geodesics $\pg_i$ in the covering above are maximal with respect to inclusion.

Indeed, we begin by the sub-geodesics containing the vertex $\pg_-$. Consider two such sub-geodesics. If one is contained in the other, by the $C'(\lambda)$--condition and the fact that $\eta \geqslant 2\lambda $ it follows that both are subwords of the same relator $r_1\in R.$
 Therefore we take the longer of the two subwords and we select it as the first term $\pg_1$ of the new covering. The endpoint $(\pg_1)_ +$ must be contained in another sub-geodesic $\pg_u$. The sub-geodesic $\pg_1 \sqcup \pg_u$ cannot be labeled by a word in $ S^\eta (R)$, because this would contradict the maximality of $\pg_1$. We consider the maximal sub-geodesic $\pg_2$ labeled by a word in $S^\eta (R)$ and containing $\pg_u$. Continuing this argument, we obtain a cover $\bigcup_{i\in S_1} \pg_i$ of $\pg$  for some $S_1 \subseteq S_0$  such that $\pg_i$ are maximal sub-geodesics labeled by words in $S^\eta (R)$.

For an arbitrary small $\varepsilon>0$ we have that $\pg \subseteq \bigcup_{i\in S_1} \pg_i^\varepsilon$, where $\pg_i^\varepsilon$ denotes the $\varepsilon$--neighborhood of $\pg_i$ in $\pg \, $.
Since $\pg$ has topological dimension one, there exists $S_2 \subseteq S_1$ such that $\pg \subseteq \bigcup_{i\in S_2} \pg_i^\varepsilon$ and every point in $\pg$ is contained in at most two sets $\pg_i^\varepsilon$ with $i\in S_2$.

If an edge $e$ in $\pg$ is not contained in $\bigcup_{i\in S_2} \pg_i$ then for $\varepsilon < \frac 12$ this contradicts the fact that $\{ \pg_i^\varepsilon \mid i\in S_2 \}$ cover $e$.
If a vertex in $\pg$ is not contained in $\bigcup_{i\in S_2} \pg_i,$ then the edges adjacent to it are not contained in  $\bigcup_{i\in S_2} \pg_i$ and we use the above.

We thus obtain that $\pg \subseteq \bigcup_{i\in S_2} \pg_i$ and every point in $\pg$ is contained in at most two sets $\pg_i$ with $i\in S_2$.

Assume that there exist two sequences $x_0=a, x_1, y_0, x_2, y_1,\ldots, x_{k+1}, y_k , y_{k+1}=b$ and  $x_0'=a, x_1', y_0', x_2', y_1',\ldots, x_{m+1}', y_m' , y_{m+1}'=b$, and let $k\leqslant m$. We prove by induction on $0\lqq i \lqq k+1$ that $[x_i, y_i] = [x_i', y_i']$.

First, consider case $i=0$. Then either $[x_0,y_0]\subseteq [x_0', y_0']$ or $[x_0',y_0']\subseteq [x_0, y_0]$. The assumption $\eta \geqslant2\lambda$ implies that both the label of $[x_0,y_0]$ and that of $[x_0', y_0']$ are subwords of the same relator $r$. The maximality condition implies that $[x_0,y_0]= [x_0', y_0']$.

We now assume that for some $j\gq 0$ we have  $[x_i, y_i] = [x_i', y_i']$ for  $0\lqq i \lqq j$. We have that either $[y_j, y_{j+1}] \subseteq [y_j, y_{j+1}']$ or $[y_j, y_{j+1}'] \subseteq [y_j, y_{j+1}]$. By maximality and Lemma \ref{lem:geod}, the contour $t_j$ containing the geodesic $[x_j,y_j]$ is distinct from the contour $t_{j+1}$ containing the geodesic $[x_{j+1},y_{j+1}]$, respectively  the contour $t_{j+1}'$ containing the geodesic $[x_{j+1}',y_{j+1}']$. It follows that
$$
\dist (x_{j+1}, y_j)< \lambda |r_{j+1}| \lqq \frac{\lambda }{\eta } \dist (x_{j+1}, y_{j+1})\, ,
$$
whence
$$
\dist (y_j ,y_{j+1}) > \left( 1- \frac{\lambda }{\eta } \right) \dist (x_{j+1}, y_{j+1})\gq \eta \left( 1- \frac{\lambda }{\eta } \right) |r_{j+1}|\, .
$$

Similarly, we obtain that
$$
\dist (y_j, y_{j+1}') > \eta \left( 1- \frac{\lambda }{\eta } \right) |r_{j+1}'|\, .
$$

The hypothesis $\eta \gq 2 \lambda $ implies that $\eta \left( 1- \frac{\lambda }{\eta } \right) \gq \lambda$, therefore the inclusions
$$
[y_j, y_{j+1}] \subseteq [y_j, y_{j+1}']\mbox{ or }[y_j, y_{j+1}'] \subseteq [y_j, y_{j+1}]
$$ imply that $t_{j+1} = t_{j+1}'$. The maximality of the sub-geodesics $[x_{j+1}, y_{j+1}]$ and $[x_{j+1}', y_{j+1}']$, and Lemma \ref{lem:geod} allow to conclude that $[x_{j+1}, y_{j+1}]=[x_{j+1}', y_{j+1}']$.

\endproof

\begin{convention}\label{conv:eta}
For the rest of this section, let $\eta $ be a fixed constant such that $\frac{1}{2}-2\lambda \gq \eta \gq 2\lambda$ and $\eta':=\eta-\lambda.$
\end{convention}

\begin{defn}[$\eta$--succession]\label{def:etasuccess}
We say that a sequence of contours $t_0,t_1,\dots,t_{k+1}$ is an $\eta$--\emph{succession of contours} if, for every $i$, one of the endpoints of $t_{i-1}\cap t_i$ is at distance $\gq \eta  |t_i|$ from at least one of the endpoints of $t_i \cap t_{i+1}$ (distance measured in $t_i$).
\end{defn}

\begin{cor}\label{cor:seqr}
Let $\pg$ be an $\eta$--relator-tied geodesic. Then there exists a unique $\eta$--succession of contours $t_0,t_1,\dots,t_{k+1}$ such that for the decomposition described in Lemma \ref{lem:rt1} the sub-geodesic with endpoints $x_i,y_i$ is contained in $t_i$.
\end{cor}

\begin{lem}\label{lem:rt2}
Let $a,b$ be two vertices joined by an $\eta$--relator-tied geodesic $\pg.$ Then every geodesic $\pg'$ with endpoints $a,b$ is $\eta'$--relator-tied, for $\eta' := \eta - \lambda$; moreover, $\pg'$ is contained in the  $\eta$--succession of contours $t_0,t_1,\ldots,t_{k+1}$ determined by $\pg$ according to Corollary \ref{cor:seqr}.
\end{lem}

\proof There exist successive points in the intersection $\pg \cap \pg'$,
$$z_0=a,z_1,\ldots, z_{2m-1}, z_{2m}, z_{2m+1}=b$$ such that $z_{2i}, z_{2i+1}$ are the endpoints of a connected component of $\pg \cap \pg'$, while $z_{2i+1}, z_{2i+2}$ are the endpoints of two sub-geodesics of  $\pg$ respectively $ \pg'$, composing a simple bigon.

Let $x_0=a, x_1, y_0, x_2, y_1,\ldots, x_{k+1}, y_k , y_{k+1}=b$ be the unique sequence of points on $\pg$ provided by Lemma \ref{lem:rt1}. For every $0\lqq i\lqq m-1$ consider the endpoints $z_{2i+1}, z_{2i+2}$ of a simple bigon. According to Theorem \ref{thm:classif}, \eqref{bigon}, the corresponding bigon is as in Figure \ref{fig:bigon}. Consider $\tau$, one of the contours appearing in this bigon; let $\alpha, \beta $ be the endpoints of the intersection of $\pg$ with $\tau$. By the small cancelation condition and the fact that the two sides of the bigon are geodesics it follows that the label of the sub-geodesic of $\pg$ limited by $\alpha, \beta$ is a sub-word of length $>\left(\frac{1}{2}-2\lambda\right)|\tau|$.

Let $m$ be the midpoint of the sub-geodesic of $\pg$ limited by $\alpha, \beta$. Then there exist $x_j, y_j$ separated by $m$. If the contour $\tau$ is distinct from the contour $t_j$ then $\frac{1}{2} \left( \frac{1}{2}-2\lambda \right) |\tau | <  \frac{1}{2}\dist (\alpha , \beta )<  \lambda |\tau |$, whence $\lambda >  \frac 18$, a contradiction. It follows that $\tau = t_j$, hence $\alpha =x_j$ and $\beta = y_j$. Thus, the endpoints of intersections of contours of the bigon with $\pg$ compose a subsequence of the sequence $x_0=a, x_1, y_0, x_2, y_1,\ldots, x_{k+1}, y_k , y_{k+1}=b$, with the property that $x_{i+1}=y_i$.

Let $z_{2i+1} $ be an endpoint of a bigon. According to the above $z_{2i+1}$ equals some $x_j$ such that $t_j$ is the first contour in the bigon. Consider now $x_{j-1}, y_{j-1}$ and the contour $t_{j-1}\neq t_j$. Then $\dist (x_j , y_{j-1} )< \lambda |t_{j-1}|$, whence $\dist (x_{j-1}, x_j) = \dist (x_{j-1}, y_{j-1})-\dist (x_{j}, y_{j-1})>\eta |t_{j-1}| -\lambda |t_{j-1}|= \eta' |t_{j-1}|$.

We thus  found that the sub-geodesic with endpoints $z_{2i}, z_{2i+1}$ common to $\pg$ and $\pg' $ is $\eta'$--relator-tied. A sub-geodesic of $\pg'$ composing one of the simple bigons is easily seen to be $\eta'$--relator-tied as $\eta' \lqq \frac{1}{2}-2\lambda $, hence the entire of $\pg'$ is $\eta'$--relator-tied.

The fact that $\pg'$ is contained in the $\eta$--succession of contours $t_0, t_1,\dots,t_{k+1}$ is immediate from the argument above: the sub-arcs of $\pg' $ with endpoints $z_{2i}, z_{2i+1}$ are contained in $\pg$, while the sub-arcs with endpoints $z_{2i+1}, z_{2i+2}$ are covered by contours $\tau$ which are in the set $\{t_0,t_1, \dots , t_{k+1} \}$.\endproof

The goal of the following two statements is to prepare the ground for the definition of the $\eta$--\dec decomposition for a pair of vertices $a,b$.

\begin{lem}\label{cor:endpoints}
Let $a$ and $b$ be two arbitrary vertices. The endpoints of an $\eta$--relator-tied component in a geodesic joining $a,b$ are contained in any other geodesic joining $a,b$.
\end{lem}

\proof Let $\pg , \pg'$ be two geodesics with endpoints $a,b$ and let $x,y$ be the endpoints of an $\eta$--relator-tied component on $\pg$. Assume that $x$ is not on $\pg'$. Then $x$ is in the interior of one of the sides of a bigon composed by $\pg$ and $\pg'$. On the other hand, this side is $\left( \frac{1}{2}-2\lambda \right)$--relator-tied, hence the component of $\pg$ between $x,y$ is not a maximal $\eta$--relator-tied sub-geodesic, a contradiction.

It follows that $x\in \pg'$ and a similar argument shows that $y\in \pg'$. \endproof

\begin{defn}[Geodesic sequences]\label{def:between}
\begin{enumerate}
  \item We say that a vertex $p$ is \emph{between two vertices $a$ and $b$} if $\dist(a,p)+ \dist (p,b)=\dist (a,b)$. We do not exclude that $p=a$ or $p=b$.

      \me

  \item We call \emph{geodesic sequence} a finite sequence of vertices $p_1,\ldots,p_m$ such that for every $1\lqq i\lqq j\lqq k \lqq m$, $p_j$ is between $p_i$ and $p_k$.

      \me

      \item If $a,b,c,d$ is a geodesic sequence then we write $(b,c) \Subset (a,d)$ and we say that the pairs $(b,c)$ and $(a,d)$ are \emph{nested}.
\end{enumerate}

\end{defn}

\begin{lem}\label{lem:concat}
Let $p,a,q,b$ be a geodesic sequence such that $p,q$ and respectively $a,b$ are the endpoints of $\eta$--relator-tied geodesics. Then there exists an $\eta$--succession of contours that contains every geodesic joining $p$ and $b$.
\end{lem}

\proof  We denote by $[p,q]$ and respectively $[a,b]$ the $\eta$--relator-tied geodesics. Consider two arbitrary geodesics $[p,a]$ and $[q,b]$ (not necessarily contained in $[p,q]$ and respectively $[a,b]$).

In the geodesic $[p,a] \cup [a,b]$, the sub-geodesic $[a,b]$ is contained in a maximal $\eta$--relator-tied component $[a', b]$. Lemma \ref{cor:endpoints} applied to $p,b$ and the geodesic joining them $[p,q] \cup [q,b]$ implies that $a' \in [p,q]$, moreover $a'$ is on every geodesic joining $p,b$. Thus, by eventually replacing $a$ with $a'$ we may assume that $a$ is contained in every geodesic with endpoints $p,b$, in particular that $a\in [p,q]$. A similar argument allows to state that without loss of generality we may assume that $q$ is contained in every geodesic joining $p,b$, in particular $q\in [a,b]$.

By Corollary \ref{cor:seqr}, there exist two $\eta$--successions of contours,
$$t_0, t_1, \dots , t_{k+1}\; \mbox{ and }\; \tau_0, \tau_1, \dots, \tau_{m+1}
$$ such that every geodesic joining $p,q$ is contained in $\bigcup_{i=0}^{k+1} t_i$, and  every geodesic joining $a,b$ is contained in $\bigcup_{j=0}^{m+1} \tau_j$.

Consider $i$ maximal such that $a\in t_i$.

Assume $i\neq k+1$. If $\tau_0 \neq t_i$ then $[a,b]\cap \tau_0$ intersects $t_i$ in a sub-geodesic of length $<\lambda |\tau_0 |$, consequently it intersects $t_{i+1}$ in a sub-geodesic of length either at least $\lambda |\tau_0 |$ or at least $(\eta  -\lambda )|t_{i+1}|$. In both cases it follows $\tau_0 =t_{i+1}\, $, whence $a\in t_{i+1}\, ,$ which contradicts the choice of $i$.

Thus, in this case, it follows that  $\tau_0 =t_i$.

Let $\ell\gq 0$ be maximal such that $\tau_r = t_{i+r}$ for $0\lqq r\lqq \ell$. It is immediate from the definition of an $\eta$--succession that the sequence
$$
t_0,\dots, t_i=\tau_0,\dots, t_{i+\ell} = \tau_\ell, \tau_{\ell+1},\dots , \tau_{m+1}
$$
is such a succession.

An arbitrary geodesic joining $p$ and $b$ must contain $a$ and $q$, the sub-geodesic from $p$ to $a$ must be contained in $\bigcup_{j=0}^i t_j$, while the sub-geodesic from $a$ to $b$ must be contained in $\bigcup_{r=0}^{m+1} \tau_r\, $.

Assume now that $i=k+1$. Every geodesic joining $a,q$ must be contained in $t_{k+1}$.

Suppose moreover that $\tau_0\neq t_{k+1}$. Then $\dist (a,q) < \lambda \min \{ |t_{k+1}|, |\tau_0| \}$. Since the distance from $q$ to one of the endpoints of $t_k \cap t_{k+1}$ is at least $\eta |t_{k+1}|$, the same is true about one of the endpoints of $\tau_0 \cap \tau_1$, since it will be situated after $q$ on a geodesic from $a$ to $b$. Therefore, in this case
$$t_0, t_1, \dots , t_{k+1}, \tau_0, \tau_1, \dots, \tau_{m+1}
$$ is an $\eta$--succession of contours.

Given an arbitrary geodesic joining $p$ and $b$, the sub-geodesic from $p$ to $q$ is in $\bigcup_{j=0}^{k+1} t_j$, the sub-geodesic from $a$ to $b$ is in $\bigcup_{r=0}^{m+1} \tau_r\, $.

Suppose that $\tau_0= t_{k+1}$. As before, the fact that $q$ is at distance $\gq \eta |\tau_0|$ from one of the endpoints of $t_k \cap t_{k+1}$ implies that one of the endpoints of $\tau_0 \cap \tau_1$ satisfies the same. Therefore, $$t_0, t_1, \dots , t_{k+1}= \tau_0, \tau_1, \dots, \tau_{m+1}
$$ is an $\eta$--succession of contours, and an argument as above shows that it contains every geodesic joining $p$ and $b$.\endproof

\comment
\begin{cor}\label{cor:seqrt}
If $a$ and $b$ are two points joined by an $\eta$--relator-tied geodesic, where $\eta  \gq 2\lambda\, $, then there exists an $\eta$--succession of contours $t_0,t_1,\ldots,t_{k+1}$ such that every geodesic with endpoints $a,b$ is covered by $t_0,t_1,\ldots,t_{k+1}$.
\end{cor}
\endcomment

\begin{rmk}\label{rem:concatN}
The statement of Lemma \ref{lem:concat} can be generalized as follows: if
$$p_0, p_1, q_0, p_2, q_1,\ldots, p_{k+1}, q_k , q_{k+1}$$ is a geodesic sequence such that $p_i,q_i$ are the endpoints of $\eta$--relator tied geodesics for $i\in \{ 0,1,\dots, k+1\}$, then there exists an $\eta$--succession of contours containing every geodesic from $p_0$ to $q_{k+1}$.

The proof adapts the argument of Lemma  \ref{lem:concat}, and we leave it as an exercise to the reader.
\end{rmk}

\begin{thm}\label{lem:intermediate}
For every pair of vertices $a,b$ in the Cayley graph  of $G$ there exists a finite geodesic sequence
$$
z_0=a,y_1, z_1, y_2, z_2,\ldots,y_m, z_m,b=y_{m+1}\, ,
$$
a sequence of $\eta$--compulsory geodesics $[z_0,y_1], [z_1, y_2],\ldots,[z_i,y_{i+1}],\ldots,[z_m, y_{m+1}]$ and a sequence of $\eta$--successions of contours
$$
t_1^{(i)},\ldots,t_{k_i}^{(i)},\, i\in \{1,2,\ldots,m \}
$$
such that $y_i \in t_1^{(i)}$, $z_i \in t_{k_i}^{(i)}$ and every geodesic joining $a,b$ is contained in
\begin{equation}\label{eq:gab}
[a,y_1]\cup \bigcup_{j=1}^{k_1} t_j^{(1)}\cup [z_1, y_2]\cup \bigcup_{j=1}^{k_2} t_j^{(2)} \cup \cdots \cup [z_{i-1},y_{i}]\cup \bigcup_{j=1}^{k_i} t_j^{(i)}\cup [z_i,y_{i+1}]\cup \cdots \cup [z_m, b]\, .
\end{equation}
\end{thm}

\begin{defn}[$\eta$--\dec decomposition]
We say that the sequence
\begin{equation}\label{eq:stdec}
(a,y_1),\backslash y_1,z_1/, (z_1,y_2), \backslash y_2, z_2/, \ldots,\backslash y_m, z_m/, (z_m,b)
\end{equation}
is \emph{the $\eta$--\dec decomposition for the pair} $a,b$.
\end{defn}

\begin{figure}[htb]
\centering
\includegraphics[scale=0.8]{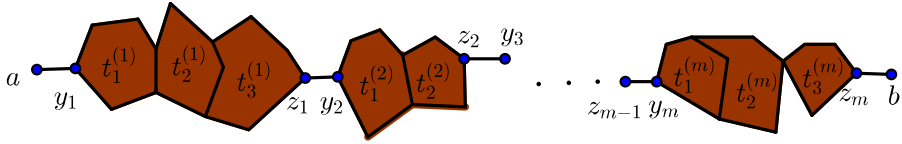}
\caption{The $\eta$--\dec decomposition for the pair $a,b$.}
\label{fig:Gab}
\end{figure}

\begin{notation}\label{notat:Gab}
For an arbitrary pair of vertices $a,b$ we denote by $\G^\eta (a,b)$ the set described in Theorem \ref{lem:intermediate},  see \eqref{eq:gab} and Figure \ref{fig:Gab}.
\end{notation}

\noindent \textit{Proof of Theorem \ref{lem:intermediate}.}   If $a,b$ is an $\eta$--compulsory pair then there is nothing to prove. Assume therefore that there exists a geodesic joining $a,b$ with an $\eta$--relator-tied component.
Let $p_1, q_1,\ldots,p_h, q_h$ be all the pairs of points that appear as endpoints of $\eta$--relator-tied components in some geodesic joining $a,b$.
Let $\pg$ be an arbitrary geodesic joining $a,b$. According to Lemma \ref{cor:endpoints}, $\pg$ contains all points $p_1, q_1,\ldots,p_h, q_h$.
 The order in which these points appear is independent of the choice of $\pg$, since it is only determined by metric relations.

 We consider the union $\bigcup_{i=1}^h [p_i, q_i]$, where $[p_i, q_i]$ denotes here the sub-geodesic of $\pg$ with endpoints $p_i, q_i$. The connected components of this union are sub-geodesics $[y_1,z_1],\ldots,[y_m, z_m]$ appearing on $\pg$ in this order. Note that $y_i\in \{ p_1,\ldots,p_h\}$ and that $z_i\in \{ q_1,\ldots,q_h\}$. In particular, both the points and the order are independent of the choice of the geodesic $\pg$.

It remains to apply Lemma \ref{lem:concat} and Remark \ref{rem:concatN}.\hspace*{\fill} $\Box$

\begin{cor}\label{cor:intermediate}
For every pair of points $a,b$ at distance $d>0$ and every $x\lqq d$ there exist at most $2$ points $p$ with the property that $a,p,b$ is a geodesic sequence and $\dist (a,p)=x$.
\end{cor}

\begin{figure}[htb]
\centering
\includegraphics[scale=0.8]{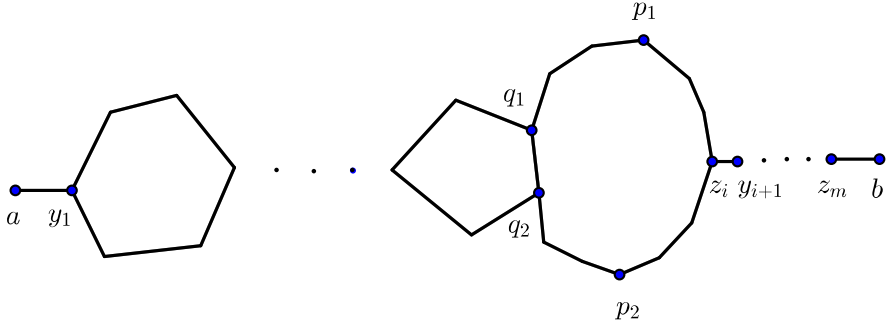}
\caption{Example of pair $a,b$ with two points between them at distance $x$ from $a$.}
\label{fig:IntermediateP}
\end{figure}

See Figure \ref{fig:IntermediateP} for an example where there exist two points $q_1, q_2$ between $a$ and $b$, at distance $x-3$ from $a$, and two points $p_1, p_2$ between $a$ and $b$, at distance $x$ from $a$.

\me
\begin{rmk}\label{rem:stdec}
Note that
\begin{enumerate}
  \item\label{rs1} according to the above every geodesic with endpoints $y_m, z_m$ is $\eta'$--relator-tied, in particular, it is non-trivial;
  \item\label{rs2} due to the maximality condition defining the pairs $y_i, z_i$, we have that $z_i \neq y_{i+1}$ for every $1\lqq i \lqq m-1$;
  \item\label{rs3} on the other hand, in the pairs $(a,y_1), (z_m,b)$ the endpoints may coincide.
\end{enumerate}
\end{rmk}

\begin{rmk}\label{rem:restr}
If $(p,q)\Subset (x,y) \Subset (a,b)$ with $p,q,x,y$ points in $\{p_1, q_1,\ldots,p_k, q_k \}$, and if $p,q$ are endpoints of an $\eta$--relator-tied component in a geodesic joining $a,b$, then $p,q$ are endpoints of an $\eta$--relator-tied component in a geodesic joining $x,y$.

This simply follows from the fact that a geodesic $\pg$  joining $a,b$ and on which $p,q$ bound an $\eta$--relator-tied component must also contain $x,y$, see Lemma \ref{cor:endpoints}.
\end{rmk}

\begin{defn}[$\eta$--relator covered pair]
If the $\eta$--\dec decomposition of a pair $a,b$ is $\backslash a,b /\, ,$ then we call such a pair an $\eta$--\emph{relator covered pair}.
\end{defn}

\comment

\begin{lem}\label{lem:compulsory}
Consider an arbitrary pair $z_i, y_{i+1}$ with $i\in \{1,2,\ldots,m-1 \}$.
\begin{enumerate}
  \item\label{comp1} There exists a unique geodesic $[z_i, y_{i+1}]$ joining $z_i, y_{i+1}$, and this geodesic is $\eta$--compulsory.
  \item\label{comp2} Every geodesic with endpoints $a,b$ contains $[z_i, y_{i+1}]$.
\end{enumerate}

\end{lem}

\proof \eqref{comp1} If there would exist two geodesics with endpoints $z_i, y_{i+1}$ then they would form at least one non-trivial bigon, hence they would contain \ertc, a contradiction.

\eqref{comp2} Every geodesic with endpoints $a,b$ contains $z_i, y_{i+1}$ and a geodesic joining them, and we apply the uniqueness of such a geodesic from \eqref{comp1}.\endproof
\endcomment

\begin{defn}[Compulsory vertices]\label{def:cpoints}
Given an $\eta$--\dec decomposition
$$(a,y_1),\backslash y_1,z_1/, (z_1,y_2), \backslash y_2, z_2/, \ldots,\backslash y_m, z_m/, (z_m,b)$$
of a pair $a,b$, we call the vertices between $z_i, y_{i+1}$ for some $i\in \{1,2,\ldots,m-1 \}$ $\eta$--\emph{compulsory vertices}.
\end{defn}

Clearly, every geodesic with endpoints $a$ and $b$ must contain all the compulsory vertices.

\begin{defn}[Prefixes and suffixes]\label{def:betweenpref}
Given an element $h\in G$ we denote by $P(h)$ (standing for \emph{prefixes of }$h$) all the elements
between $1,h$ and by $S(h)$ (standing for \emph{suffixes of }$h$) all the elements of the form $x\iv h$ for $x\in P(h)$.
\end{defn}

Note that the two sets $P(h)$ and $S(h)$ depend on the fixed generating set $A.$

\begin{defn}[Compulsory and $\eta$--relator-covered elements]\label{def:elements}
Let $h\in G$.
\begin{itemize}
  \item If $h$ is joined to $1$ by at least one \ertgs then we call $h$ an $\eta$--\emph{relator-tied element}.
  \item If the pair $1,h$ has the $\eta$--\dec decomposition $(1,h)$ (hence, there exists only one geodesic joining $1,h$, composed of compulsory vertices), then
we call $h$ an $\eta$--\emph{compulsory element}.
  \item  If the pair $1,h$ has the $\eta$--\dec decomposition $\backslash 1,h /$ then
we call $h$ an \emph{$\eta$--relator-covered element}.
\end{itemize}
\end{defn}

\begin{notation}\label{notat:comp}
We denote by $RT^\eta$ the set of $\eta$--relator-tied elements.
We denote by $\calc^\eta$ the set of $\eta$--compulsory elements in $G$ and by $RC^\eta$ the set of $\eta$--relator-covered elements.
\end{notation}

\begin{rmk}
The fact that $h$ is $\eta$--relator-covered does not mean that there exists an \erts geodesic labeled by $h$, it only means that every geodesic $[a,b]$ labeled by $h$ contains a family of successive vertices $y_0=a, y_1, z_0, y_2,z_1,\ldots,$ $y_m, z_{m-1}, z_m =b$ such that for every $i$ there exists an \erts geodesic with endpoints $y_i,z_i$. In particular, by Lemma \ref{lem:rt2}, every geodesic labeled by $h$ is \eprt.

\end{rmk}

An algebraic version of Theorem \ref{lem:intermediate} is the following.

\begin{thm}\label{cor:sdecomp}
 Every element $g\in G$ can be written uniquely as a product
 \begin{equation}\label{eq:dec-elem}
  g= \alpha_1\beta_1\alpha_2\beta_2\dots \alpha_m \beta_m \alpha_{m+1}\, ,
 \end{equation}
 such that
\begin{itemize}
  \item $\beta_i$ are non-trivial $\eta$--relator-covered elements;
  \item $\alpha_i$ are compulsory elements (non-trivial with the possible exception of $\alpha_1,$ $\alpha_{m+1}$);
  \item the vertices $y_i= \alpha_1\beta_1\dots \alpha_i$ and $z_i = \alpha_1\beta_1\dots \alpha_i \beta_i$ compose the geodesic sequence determining the $\eta$--\dec decomposition of the pair $1,g$.
\end{itemize}
\end{thm}

\begin{notation}\label{notat:Gh}
Given an arbitrary element $h\in G$ we denote by $\G^\eta (h)$ the set $\G^\eta (1,h)$ as described in Notation \ref{notat:Gab}.
\end{notation}

\comment
We write that the $(2m+1)$--tuple $(\alpha_1,\beta_1,\alpha_2,\beta_2,\dots ,\alpha_m, \beta_m, \alpha_{m+1})$ is in $\dd_m$ if it corresponds to a decomposition as above. If $\alpha_1$ or $\alpha_{m+1}$ are trivial we omit them, but still say that the $(2m-1)$--tuple or the $(2m)$--tuple is in $\dd_m$.
\endcomment

\comment
\begin{notation}\label{notat:dmlA}
Given a finite sequence of non-negative integers $A = (a_1,b_1,\ldots,a_s)$ or $B = (a_1,b_1,\ldots,a_s, b_s)$ we denote by $\Lambda (A)$, respectively $\Lambda (B)$ the set of elements $g$ in $G$ with a decomposition as in \eqref{eq:dec-elem} such that $m\gq s$ and $|\alpha_i| =a_i$ for all $1\lqq i\lqq s$, $|\beta_i| =b_i$ for all $1\lqq i\lqq s-1$ (respectively,  all $1\lqq i\lqq s$).

When $m=s$ and $A$, respectively $B$, lists the lengths of all the components involved in the decomposition of $g$, we use the notation $\Lambda_t (A)$, respectively $\Lambda_t (B)$.
\end{notation}

\endcomment

\begin{notation}\label{notat:dd_i}
Given $i\in \N,\, i\gq 2$, we denote by $\dd_i$ the set of $i$--tuples
$$(a_1,a_2,\ldots,a_{i-1}, b)$$ such that for the element $g= a_1a_2\cdots a_{i-1} b$ the elements $a_1,a_2,\ldots,a_{i-1}$ are the first $i-1$ elements in the \dec decomposition of $g$ as described in Theorem \ref{cor:sdecomp}.
\end{notation}

\comment
\begin{cor}\label{cor:prefixrc}
Let $\beta$ be a non-trivial $\eta$--relator-covered element. Then every non-trivial element $\sigma \in P(\beta )$ is either compulsory or it has a decomposition as above of the form
$$
\sigma = \beta_0 \alpha_0\, ,
$$ where $\beta_0 \in P(\beta )$ is an $\eta$--relator-covered element and $\alpha_0$ is compulsory (possibly trivial).
\end{cor}

\proof By Lemma \ref{cor:endpoints}, the pairs with endpoints of $\eta$--relator-tied components in some geodesic joining $1$ and $\beta$ appear also on every other geodesic joining $1$ and $\beta$. We consider the list $p_1, q_1,\ldots,p_k, q_k$ of all such pairs of vertices, except for the ``nested'' pairs, that is the pairs with endpoints $p,q$ for which there exists another such pair $p',q'$ such that $1, p', p, q, q' , \beta$ is a geodesic sequence.

Consider the vertex in $\{ p_1, q_1,\ldots,p_k, q_k \}$ that is between $1$ and $\sigma $ and is nearest to $\sigma$.

\textit{Case 1.} Assume that this point is a $q_i$.  If $q_i, \sigma$ is an $\eta$--compulsory pair then we are done.

Assume that $q_i,\sigma $ is not an $\eta$--compulsory pair, and consider a geodesic $\pg$ joining $1$ and $\beta$ and containing $q_i$ and $\sigma$. An \ertgs subsegment of $\pg$ between $q,\sigma $ must extend to a maximal \ertgs contained in $\pg$ with start before $q$ and ending after $\sigma$. Otherwise the choice of $q_i$ would be contradicted. Let $[a,b] \subset \pg$ be this maximal \ertg.

Consider the unique sequence of points in $[a,b]$ defined by Lemma \ref{lem:rt1}
$$x_0=a, x_1, y_0, x_2, y_1,\ldots, x_{m+1}, y_m , y_{m+1}=b\, .$$

\me

\textit{Case 2.} Assume now that the nearest point to $\sigma$ between it and $1$ is a $p_i$.

According to the hypothesis that $1, \beta$ is an $\eta$-relator-covered element, and the hypothesis of non-nesting it follows that either $p_i$ equals some $q_j$, ir that there exists a pair $p_j,q_j$ with $p_j$ between $1,p_i$ and $q_j$ between $\sigma$ and $q_i$. In the first case we argue as in Case 1. Assume therefore that we are in the second case. An argument as in Case 1 with $a,b$ replaced by $p_j,q_j$ \ldots

\endproof
\endcomment

The following lemma will be crucial for the results in Section \ref{sec:qm} on quasi-homo\-mor\-phisms.

\begin{lem}\label{lema:T1ab}
Let $\lambda \leqslant \frac{1}{10}$ and let $\eta \gq 3\lambda$.

Every $\eta$--succession of contours $t_0,t_1,\dots,t_{k+1}$ is totally geodesic: if $a,b$ are two vertices in $\bigcup_{i=0}^{k+1} t_i$ then every geodesic joining $a$ and $b$ is contained in $\bigcup_{i=0}^{k+1} t_i$.
\end{lem}

\proof Without loss of generality we assume that $a\in t_1\setminus t_2$ and that $b\in t_{k+1}\setminus t_k$. Otherwise, assuming that $a$ appears before $b$ in the succession, we consider the largest $i$ such that $t_i$ contains $a$ and the smallest $j$ such that $t_j$ contains $b$ and take the succession $t_i,t_{i+1},\dots, t_{j-1}, t_j$ instead of the initial one.

Let $\pg $ be a geodesic joining $a$ and $b$. We argue for a contradiction and assume that $\pg$ is not contained in $\bigcup_{i=0}^{k+1} t_i$. Without loss of generality we assume that $\pg$ intersects  $\bigcup_{i=0}^{k+1} t_i$ only in its endpoints (otherwise, we replace $\pg$ by a sub-geodesic with this property).

Let $\pgot$ be a topological arc joining $a$ and $b$ in $\bigcup_{i=0}^{k+1} t_i$ and of minimal length. By the Greendlinger Lemma, there exists a contour $\tau$ such that one of the connected components of its intersection with $\pgot \cup \pg$ has length $>(1-3\lambda )|\tau |$. If $\tau = t_i$ for some $i$ then by the hypothesis on $\pg$, $\tau$ intersects $\pgot$ in a connected component of length $>(1-3\lambda )|\tau |$. Then $\pgot$ can be shortened by a length of $(1-6\lambda )|\tau |$, which contradicts the choice of $\pgot$ as an arc of minimal length joining $a$ and $b$ in $\bigcup_{i=0}^{k+1} t_i$.

We therefore assume that $\tau \not\in \{t_0,t_1,\dots,t_{k+1}\}$. Since $\pg$ is a geodesic, it follows that $\tau$ intersects $\pgot$ in a subarc of length $>\left( \frac{1}{2}-3\lambda \right)|\tau |$.

On the other hand, $\pgot$ contains a succession of vertices $$x_0=a, x_1, y_0, x_2, y_1,\ldots, x_{k+1}, y_k , y_{k+1}=b$$
such that the sub-arcs with endpoints $x_i, y_i$ with $i\in \{ 0,1,\ldots,k+1\}$ are labeled by words in $S(R)$, which are moreover in $S^\lambda (R)$ if $i\neq 0, k+1$. Therefore the connected component of the intersection $\tau \cap \pgot$ cannot contain a pair $x_i,y_i$ with $i\in \{1,\ldots,k\}$. It follows that it can intersect at most two consecutive sub-arcs with endpoints $x_i, y_i$ with $i\in \{ 0,1,\ldots,k+1\}$, hence it is of length $<2\lambda |\tau|$. We thus obtain that $\frac{1}{2}-3\lambda < 2\lambda$, whence $\lambda >\frac{1}{10}$, a contradiction.\endproof

The following results are not used in an essential manner in our arguments, but they complete nicely the description of geodesics in small cancelation groups.

\begin{lem}\label{sub-ert}
Let $\pg$ be an $\eta$--relator-tied geodesic and let $\p$ be a sub-geodesic in it. Then $\p$ is either an $\eta$--compulsory geodesic, or it is the concatenation of three sub-geodesics $\p = \p_c \sqcup \p_0 \sqcup \p_c'$, where $\p_c ,\p_c'$ are $\eta$--compulsory and contained in a contour (possibly either one of them or both trivial) and $\p_0$ is an $\eta$--relator-tied component of $\p$.
\end{lem}

\proof

Let $a,b$ be the endpoints of $\pg$. With the previous convention $\pg = [a,b]\, $.

\me

\emph{Step 1.} Let us first assume that $\p =[a, \sigma ]$, with $a, \sigma , b$ a geodesic sequence.

Let $x_0=a, x_1, y_0, x_2, y_1,\ldots, x_{k+1}, y_k , y_{k+1}=b$ be the unique sequence of points on $\pg$ defined by Lemma \ref{lem:rt1}.

Assume that $\sigma $ is in between a pair $y_j, x_{j+2}$. If the word labeling the geodesic $[x_{j+1}, \sigma ]$ is contained in $S^\eta (R)$ then $\p$ is $\eta$--relator-tied.

If the word labeling $[x_{j+1}, \sigma ]$ is not in  $S^\eta (R)$ (while it is still a sub-word of the relator labeling the contour $t_{j+1}$) then the pair $x_{j+1}, \sigma$ is $\eta$--compulsory. This implies that the required decomposition is
 $\p = [a, y_j]\sqcup [y_j, \sigma ]\, $.

Assume now that $\sigma $ is  in between a pair $x_{j+1}, y_j$. If $[x_j, \sigma]$ is labeled by a word in $S^\eta (R)$ then $\p$ is $\eta$--relator-tied;
while in the opposite case the geodesic $[x_j, \sigma]$ is $\eta$--compulsory, and the conclusion holds with the decomposition $\p =[a ,y_{j-1}] \sqcup [y_{j-1} , \sigma ]$.

\me

\emph{Step 2.} Assume now that $\p = [\varrho , \sigma ]$, where $a, \varrho , \sigma ,b$ is a geodesic sequence. According to Step 1, $[a, \sigma ] = [a, \mu ] \sqcup [\mu , \sigma ]\, $, where $[a, \mu ]$ is an $\eta$--relator-tied component and $[\mu , \sigma ]$  is  $\eta$--compulsory (possibly trivial) and contained in a contour. If $\varrho \in [\mu , \sigma ]$ then $\p$ is  $\eta$--compulsory. If $\varrho \in [a, \mu ]$ then by reversing the order on $[a, \mu ]$ and applying Step 1 we obtain that $[\varrho , \mu ] = [\varrho , \nu ] \sqcup [\nu , \mu ]\, ,$ where $[\varrho , \nu ]$ is $\eta$--compulsory (possibly trivial) and contained in a contour, and $[\nu , \mu ]$ is an $\eta$--relator-tied component. It follows that
$$\p = [\varrho , \mu ] \sqcup [\mu , \sigma ] =  [\varrho , \nu ] \sqcup [\nu , \mu ]\sqcup [\mu , \sigma ]\,
$$ is the required decomposition.
\endproof

 \begin{lem}
For each pair $\backslash y_j, z_j/$ in an $\eta$--\dec decomposition there exists a geodesic sequence
\begin{equation}\label{eq:seqyz}
 p_1'= y_j, p_2', q_1', p_3', q_2',\ldots,p_n', q_{n-1}', q_n'=z_j\, ,\mbox{ for some } n=n(j)\, ,
\end{equation}

such that:
\begin{itemize}
  \item $(p_s', q_s')$ are maximal with respect to the partial order relation $\Subset$;
  \item $p_{\ell +1}',  q_\ell'$ bound $\eta$--relator-tied sub-geodesics both in the $\eta$--relator-tied geodesic joining $p_\ell', q_\ell'$ and in the $\eta$--relator-tied geodesic joining $p_{\ell +1}', q_{\ell +1}'$.
\end{itemize}
 \end{lem}

 \proof By definition, $[y_j, z_j] =\bigcup_{i\in I_j} [p_i, q_i]\, $. Without loss of generality we assume that there are no nested pairs among the $(p_i, q_i)$ with $i\in I_j$, in other words each pair $(p_i, q_i)$ is maximal with respect to the partial order relation $\Subset$. Proceeding as in the proof of Lemma \ref{lem:rt1} we  also assume that, after selecting a subset in $I_j$, every point on a (every) geodesic joining $y_j, z_j$ is between at most two pairs $(p_i, q_i)$. It then follows that the set of pairs indexed by $I_j$ compose a geodesic sequence as in \eqref{eq:seqyz}. We set $p_i' := p_i$ and $q_i' :=q_i$.

Consider now two consequent pairs that overlap: two pairs $(p_i', q_i' )$ and $(p_{i+1}' , q_{i+1}')$ such that $p_i', p_{i+1}' , q_i' , q_{i+1}'$ is a geodesic sequence.

 By definition, there exists a geodesic $\pg$ joining $a,b$ such that $p_i', q_i'$ are the endpoints on it of an $\eta$--relator-tied component.  Given two points $x,y\in \pg,$ we denote by $[x,y]$ the sub-geodesic of $\pg$ with endpoints $x,y$.

  We likewise know that there exists a geodesic $\p$ such that $p_{i+1}' , q_{i+1}'$ bound an \erts component on $\p$. According to the above, $\p$ must contain $q_i'$. In what follows, for $x,y$ in $\p$ we denote by $\overline{x,y}$ the sub-arc of $\p$ with endpoints $x,y$.

 We have that $p_{i+1}' \in [p_i', q_i']$. Lemma \ref{sub-ert} implies that either $[p_{i+1}' , q_i' ]$ is an $\eta$--compulsory component contained in a contour, or $[p_{i+1}' , q_i' ] = [p_{i+1}' , x] \sqcup [x, q_i']$, where $[p_{i+1}' , x] $ is an $\eta$--compulsory component contained in a contour (possibly trivial) and $[x, q_i']$ is an $\eta$--relator-tied component.

 Assume that $[p_{i+1}' , q_i']$ is an $\eta$--compulsory component contained in a contour. Then the geodesic $\p$ must also contain $[p_{i+1}' , q_i' ] \subset \pg$. By replacing on $\p$ the sub-arc with endpoints $a,p_{i+1}'$ by $[a, p_{i+1}' ]\subset \pg$ we obtain a new geodesic $\rgot$ joining $a,b$ such that $p_i'$ and $q_{i+1}'$ are the endpoints of an \erts sub-geodesic. It follows that $(p_i', q_{i+1}') \Subset (\alpha , \beta)$, where $\alpha , \beta$ are the endpoints on $\rgot$ of an \ertc. In particular $(\alpha , \beta ) = (p_\ell, q_\ell )$ for some $\ell\in I_j$, and $(p_i', q_i') \Subset (p_\ell, q_\ell )$. This contradicts the fact that we have considered pairs maximal with respect to $\Subset$.

 Assume that $[p_{i+1}' , q_i' ] = [p_{i+1}' , x] \sqcup [x, q_i']$, where $[p_{i+1}' , x] $ is an $\eta$--compulsory component contained in a contour (possibly trivial) and $[x, q_i']$ is an $\eta$--relator-tied component. Since $q_i'\in \p$ and $[x,q_i']$ is an $\eta$--relator-tied component between $p_{i+1}'$ and $q_i'$ it follows that $x\in \p$, hence $[p_{i+1}' , x] \subset \p$.
  There exists $y\in \p$ such that $\overline{p_{i+1}', y}$ is labeled by a word in $S^\eta (R)$ and it is contained in a contour $t$. If $y\in \overline{p_{i+1}' , x} = [p_{i+1}' , x]$ then the contour $t$ intersects a distinct contour in a sub-arc of length $\gq \eta |t|$, a contradiction. Hence we must have that $x\in \overline{p_{i+1}', y}$.

  According to the small cancelation condition $\overline{x,y}$ has length $>\left( 1- \frac{\lambda }{\eta }\right)$ of the length of $\overline{p_{i+1}', y}$, so at least $\eta \left( 1- \frac{\lambda }{\eta }\right) |t|$. This implies that if $\eta \left( 1- \frac{\lambda }{\eta }\right) \gq\lambda , $ equivalently $\eta \gq2 \lambda $, then by Lemma \ref{lem:rt2}, $t$ must be the first contour for the pair $x, q_i'$. But this implies that $p_{i+1}' = x$.

  Similarly, we argue that $\overline{p_{i+1}', q_i'}$ is an $\eta$--relator-tied geodesic. \endproof

\section{Proof of the Rapid Decay Property}\label{sec:main}
We apply the \dec decomposition of elements in $G$ that we have developed in Section \ref{sec:stdec} to prove
Theorem \ref{thm:main}.

\begin{cvn}\label{conv:supp-eta}
Throughout this section we fix arbitrary constants $R\gq r \gq 0, p\in [R-r , R+r]$ and  an arbitrary function $f\in \R_+G$
    with support in $S(1,r)$.

We fix $\lambda \lqq \frac{1}{10}$. Parts of the argument work for $\lambda \lqq \frac{1}{8}$ as well; the choice  $\lambda \lqq \frac{1}{10}$ is required to
conclude cases $T_2, (I)$ and $T_2, (III)$ below.

We  fix a constant $\eta$ such that $\frac{1}{2} -2 \lambda \gq \eta\gq 2\lambda$.
\end{cvn}

    Consider the operator
    $$
    T_f : l^2(S(1,R)) \to l^2(S(1,p)),\, T_f (g) = (f*g)_p\, .
    $$

Our goal, according to Lemma \ref{equiv1}, is to prove that
\begin{equation}\label{eq:operator}
\|T_f\|\lqq P(r) \|f\|\, ,
\end{equation}
where $\|T_f\|$ denotes the operator norm and $P$ is a polynomial independent of the choice of $R,r,p$ and $f$.

\begin{figure}[htb]
\centering
\includegraphics[scale=0.9]{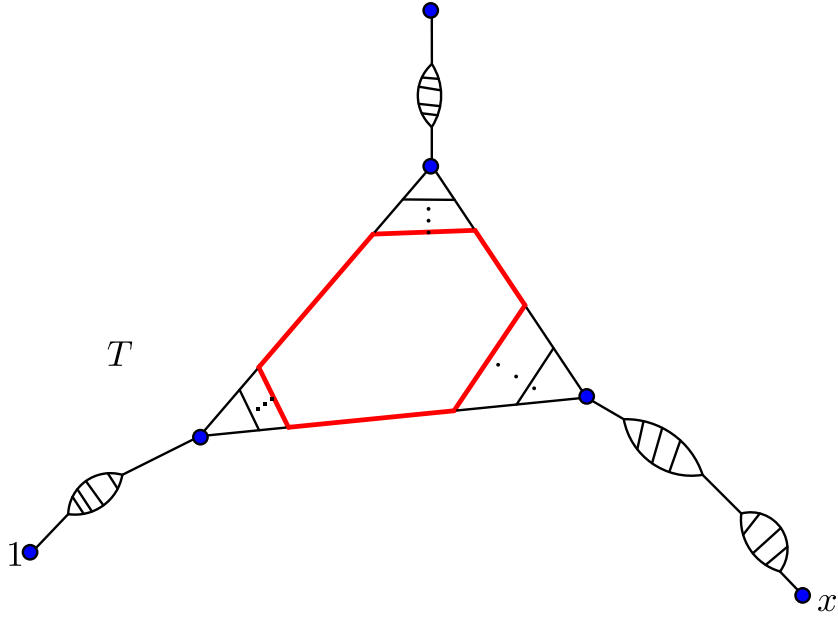}
\caption{A geodesic triangle in $G$.}
\label{fig:gentriangle}
\end{figure}

Given three elements $a,b,x$ in $G$ such that $ab=x$ we denote by $\Delta (x,a,b)$ the set of geodesic triangles in the Cayley graph with vertices $1,x,a$ (hence with edges labeled by shortest words representing $x,a,b$). See Figure \ref{fig:gentriangle} for a general picture.

Given an arbitrary geodesic triangle $\trian$ in $\Delta (x,a,b)$ with edges denoted by $[1,x]$, $[1,a]$ and $[a,x]$ (the last labeled by a shortest word representing $b$), there exists a central simple triangle defined as follows. Let $C_1$ be the farthest from $1$ point in $[1,x]\cap [1,a]$,  $C_a$ be the farthest from $a$ point in $[a,1]\cap [a,x]$, and $C_x$ be the farthest from $x$ point in $[x,1]\cap [x,a]$.

The triangle composed by $[C_1, C_x] \subseteq [1,x],\, [C_1, C_a] \subseteq [1,a],\, [C_a, C_x] \subseteq [a,x]$ is simple (possibly trivial). We call it the \emph{central triangle of} $\trian$, and we denote it by $\theta (\trian )$. If two of the vertices of $\theta (\trian )$ are equal, then by definition the third vertex is also equal to the other two, hence the central triangle is trivial.
If $\theta (\trian )$ is non-trivial then it has pairwise distinct vertices and it has one of the shapes described in Theorem \ref{thm:classif},  see Figures \ref{fig:bigon} and \ref{fig:triangle}.

Let $\Delta^i (x,a,b)$ denote the set of geodesic triangles $\trian \in \Delta (x,a,b)$ with central triangle $\theta (\trian )$ non-trivial and of the form $T_i$, where $i=1,2,3,4$ and $T_i$ is the respective shape in  Figures \ref{fig:bigon} and \ref{fig:triangle}. We denote by $\Delta^0 (x,a,b)$ the set of geodesic triangles $\trian $ in $\Delta (x,a,b)$ with central triangle $\theta (\trian )$ trivial. Clearly,
$$
\Delta (x,a,b) = \Delta^0 (x,a,b) \cup \Delta^1 (x,a,b)\cup \Delta^2 (x,a,b)\cup \Delta^3 (x,a,b)\cup \Delta^4 (x,a,b)\, ,
$$
but the union is not disjoint. We therefore have that for an arbitrary $x$ with $|x|=p$
$$
T_f(g)(x)=\sum_{ab=x} f(a)g(b) \lqq \sum_{i=0}^4\sum_{ab=x, \Delta^i (x,a,b) \neq \emptyset} f(a)g(b)\, .
$$

It suffices to prove \eqref{eq:operator} for the operators $T^i_f$, with $i=0,1,2,3,4,$ defined by
\begin{equation}\label{eq:ti}
T^i_f(g)(x)=\sum_{ab=x, \Delta^i (x,a,b) \neq \emptyset} f(a)g(b)\, .
\end{equation}

We have to thus estimate the norm of $T^i_f(g)(x)$, and prove that
\begin{equation}\label{eq:norm}
\sum_{|x|=p} \left( \sum_{ab=x, \Delta^i (x,a,b) \neq \emptyset} f(a)g(b)  \right)^2 \lqq P(r) \| f\|^2 \| g\|^2\, .
\end{equation}

For convenience, in what follows, for a given element $x\in G$ of length $p$ we denote by $S^i (x) $ the term $\sum_{ab=x, \Delta^i (x,a,b) \neq \emptyset} f(a)g(b)\, $. By abuse of notation we use the same even when we further restrict the case $\Delta^i (x,a,b) \neq \emptyset$, i.e. when we actually consider a sum with less terms.

We discuss each of the cases $\Delta^i (x,a,b) \neq \emptyset$, with $i=0,1,2,3,4$; we denote each case by $T_i$, corresponding to the type of the central triangle as represented in Figures \ref{fig:bigon} and \ref{fig:triangle}. The order in the proof will not follow the order of the indices $i$ of $T_i$, it will follow the logic of the argument.

\medskip

\noindent \textbf{Case $T_0$. $\Delta^0 (x,a,b) \neq \emptyset$}.

In other words, we work with products $ab=x$ such that there exists a point $\theta $ that is simultaneously between $(1,x)$, $(1,a)$ and $(a,x)$.

Note that necessarily, $\dist (1, \theta ) = \frac{1}{2} (r+p-R), \dist (a, \theta ) = \frac{1}{2} (r+R-p)$ and $\dist (\theta , x )= \frac{1}{2} (R+p-r)$.

The sum on the left hand side of \eqref{eq:norm} is at most
\begin{equation}\label{eq:sum1}
\sum_{|x|=p} \left(\sum_{\theta x_1=x, |\theta |= \frac{1}{2} (r+p-R), |x_1|=  \frac{1}{2} (R+p-r)}\; \;  \sum_{|\ell |= \frac{1}{2} (r+R-p)} f(\theta \ell )g(\ell\iv x_1)  \right)^2 \, .
\end{equation}

According to Corollary \ref{cor:intermediate}, a given element $x$ of length $p$ has at most $2$ distinct decompositions $\theta x_1=x$ with $|\theta |= \frac{1}{2} (r+p-R), |x_1|=  \frac{1}{2} (R+p-r)$. Thus
the particular case of the Cauchy-Schwartz inequality $(a_1 + a_2)^2 \lqq 2(a_1^2 + a_2^2)$ implies that an upper bound for \eqref{eq:sum1} is obtained by moving the first sum inside brackets in \eqref{eq:sum1}, outside the brackets, and multiplying by $2$. Now the sum outside the brackets is over all $x$ of length $p$, and over all decompositions $\theta x_1=x$ of the required type. By eventually further increasing the value we may replace this double sum with a sum over all $\theta$ and $x_1$ with $|\theta |= \frac{1}{2} (r+p-R), |x_1|=  \frac{1}{2} (R+p-r) $.

We thus obtain that an upper bound of the sum in \eqref{eq:sum1} is
$$
2\sum_{|\theta |= \frac{1}{2} (r+p-R)}\; \; \sum_{|x_1|=  \frac{1}{2} (R+p-r)} \left(\sum_{|\ell |= \frac{1}{2} (r+R-p)} f(\theta \ell )g(\ell\iv x_1)  \right)^2\, .
$$

We apply the Cauchy-Schwarz inequality in the last sum above and we obtain an upper bound of the form $2\Pi_1 \, \Pi_2$, where
$$
\Pi_1=\sum_{|\theta |= \frac{1}{2} (r+p-R), |\ell |= \frac{1}{2} (r+R-p)} [f(\theta \ell )]^2
$$ and
$$
\Pi_2= \sum_{|\ell'| = 1/2(r+R -p), |x_1|=  \frac{1}{2} (R+p-r)} [g(\ell' x_1 )]^2\, .
$$

We compare $\Pi_1$ to $\|f\|^2$, by asking, given an arbitrary element $u\in G$, how many times does $[f(u)]^2$ appear in the sum $\Pi_1$, in other words in how many different ways can $u$ be written as $u= \theta \ell \, $. By Corollary \ref{cor:intermediate} each $u$ can have at most $2$ such decompositions, thus $\Pi_1 \lqq 2\|f\|^2$.
We argue likewise that $\Pi_2 \lqq 2\|
g\|^2$.

\comment
Given an arbitrary $x$,  the pair $1,x$ has the $\eta$--\dec decomposition of the form
$$(1,y_1),\backslash y_1,z_1/, (z_1,y_2), \backslash y_2, z_2/ \ldots,\backslash y_m, z_m/, (z_m,x)\, .$$

The norm of the operator $T^0_f$ is bounded by the sum of the norms of the two operators $T^{0,r}_f$ and  $T^{0,c}_f$ corresponding to the two cases when  $\theta (x,a,b)$ is in between an \erc pair or in between an $\eta$-compulsory pair. Since the two cases are not disjoint we cannot say that $T^0_f$ is the sum of $T^{0,r}_f$ and  $T^{0,c}_f$, but only have the inequality between norms mentioned above.

\medskip

\noindent \textbf{Case $T_0$, (a). $\Delta^0 (x,a,b) \neq \emptyset$ and some $\theta (x,a,b)$ is between $y_s, z_s$, distinct from  $y_s, z_s$, for an $s\gq 1$}.

Note that in particular this implies $\dist (1,y_s)\lqq \frac{1}{2} (r+p-R)\lqq r$. Since we have a number of possibilities for $s$ that is linear in $r$ we may as above split the operator $T^{0,r}_f$ as a sum of operators, and reduce to the estimate of the norm of the operator with $s\in [0,r+1]$ fixed.

Once $s$ is fixed, the number of possible choices for the distances
\begin{equation}\label{eq:distan}
 \dist(1,y_1),\dist (y_1,z_1), \dist(z_1,y_2), \dist (y_2, z_2) \ldots,\dist(z_{s-1}, y_s)
\end{equation}
is also polynomial in $r$, therefore we may again split the operator into a sum of operators, one for each choice of values for the distances in \eqref{eq:distan}, and reduce the problem of bounding the norm for the whole operator to that of each term.

We may thus assume from now on that the first sum on the left hand side of \eqref{eq:ti} is over $x$ such that on every geodesic $[1,x]$ the point at distance $\frac{1}{2} (r+p-R)$ is inside a component $y_s, z_s$ with $s$ fixed, and the distances in \eqref{eq:distan} take fixed values $a_1, b_1,a_2,b_2,\dots, a_s$, and the second sum is over $ab=x$ such that $\Delta^0 (x,a,b)\neq \emptyset$.

We use the notation $A=(a_1, b_1,a_2,b_2,\dots, a_s)$ and $\Lambda (A), \Lambda_t (A)$ introduced in \ref{notat:dmlA}.

We then use the Remark \ref{rem:stdec}, \eqref{rs1}, the Lemma \ref{sub-ert} and the notation in \ref{notat:comp} to deduce that the elements $x$ over which is taken the first sum on the left hand side of \eqref{eq:ti} can all be decomposed as products $x= h \tau \chi \chi' \tau'k$, where $|h \tau \chi| = \frac{1}{2}(r+p-R)$, and $h\in \Lambda_t (A), \tau, \tau' \in RT^{\eta'}, \chi, \chi'\in C^{\eta'}, |k| \lqq 1/2(p+R -r)$.

Then the left-hand side of \eqref{eq:norm}, with the sums restricted as above, is bounded by
$$
\sum_{h\in \Lambda_t (A), \tau, \tau' \in RT^{\eta'}, \chi, \chi'\in C^{\eta'}, |k| \lqq 1/2(p+R -r)}\left(\sum_{|\ell| = 1/2(r+R -p)} f(h\tau \chi \ell ) g(\ell\iv \chi'\tau' k ) \right)^2\, .
$$

We apply the Cauchy-Schwarz inequality in the last sum (over $\ell$) above, and we obtain the upper bound
$\Pi_1 \, \Pi_2$, where
$$
\Pi_1=\sum_{h\in \Lambda_t (A), \tau \in RT^{\eta'}, \chi \in C^{\eta'}, |\ell| = 1/2(r+R -p)} [f(h\tau \chi \ell )]^2
$$ and
$$
\Pi_2= \sum_{|\ell'| = 1/2(r+R -p), \chi'\in C^{\eta'}, \tau' \in RT^{\eta'}, |k| \lqq 1/2(p+R -r)} [g(\ell' \chi'\tau' k )]^2\, .
$$

We compare $\Pi_1$ to $\|f\|^2$, by asking, given an arbitrary element $u\in G$, how many times does $[f(u)]^2$ appear in the sum $\Pi_1$.

If $u$ can be written as a product $h\tau \chi \ell$ then by the uniqueness of the $\eta$--\dec decomposition, the factor $h$ is uniquely determined by $u$.

If $\chi$ is non-trivial then $\tau$ is an \eprtcs of $h\iv u$, hence it is also uniquely defined. Both $h$ and $h\tau$ appear on every geodesic joining $1$ and $u$ by Lemma \ref{cor:endpoints}. The point $\theta$ is inside an $\eta'$--compulsory component starting in $h\tau$. That component too is included in every geodesic joining $1$ and $u$, in particular so is $\theta$, and it is the unique point at distance $\frac{1}{2} (r+p-R)$ from $1$. Therefore $\chi$ is also uniquely determined by $u$. Thus in this case, every $u$ can be written at most in one way as a product $h\tau \chi \ell$.

If $\chi $ is trivial then $[h, h\tau ]$ is contained in an \eprtcs of the geodesic $[h, u]$. Assume that the second endpoint of that component is some $h\widetilde{\tau}$. By Lemma \ref{cor:endpoints} the element $h\widetilde{\tau}$ is contained in every geodesic $[1,u]$, and $\eta'>3\lambda$ according to Convention \ref{conv:sup-eta}.

It follows that for every geodesic $[1,u]$, the sub-geodesic between $h$ and $h\widetilde{\tau}$ is covered by an $\eta$--succession of contours as in Lemma \ref{lem:rt2}. Hence we have at most 4 distinct choices for a point $\theta$ on a geodesic $[1,u]$ at distance $\frac{1}{2} (r+p-R)$ from $1$. Hence in this case, every $u$ can be written at most in four different ways as a product $h\tau \ell$.

We conclude that $\Pi_1 \lqq 4\|f\|^2$.

We now compare $\Pi_2$ with $\|g\|^2$. In other words, for an arbitrary element $v\in G$, we study how many times  $[g(v)]^2$ appears in the sum $\Pi_2$.

By hypothesis, the point $\ell'\chi'\tau'$ is the first of the points defining a \dec $\eta$-decomposition as in \eqref{eq:stdec} for $1,v$ to appear at a distance $\gq \frac{1}{2} (r+R-p)$ from $1$. Hence for every given $v$ there is only one possible choice for this point, hence for $k$.

If $\chi'$ is non-trivial, then $\ell'\chi'$ is the endpoint of the last \eprtcs between $1$ and $\ell'\chi'\tau' = gk\iv$. It is therefore uniquely defined, which means that $\tau'$ is also uniquely defined. The geodesic $[\ell' \, ,\,  \ell'\chi']$ is part of an $\eta'$-compulsory component of $[1, \ell'\chi']$. In particular there is a unique choice for $\chi'$, hence  in this case $v$ can have at most one decomposition as $\ell'\chi'\tau' k$.

Assume that $\chi'$ is trivial. Then a geodesic $[\ell' , \ell'\tau' ]$ is part of an \eprtcs of $v$. An argument as in Corollary \ref{cor:intermediate} shows that there are at most $2$ distinct choices for a point on a geodesic $[1,v]$ at distance $\frac{1}{2} (r+R-p)$ from $1$, hence at most $2$ possible choices for $\ell'$. It follows that in this case, every $v$ can be written in at most $2$ distinct ways as $\ell'\tau'k$.

We can therefore write that $\Pi_2 \lqq 2\|g\|^2$.

\medskip

\noindent \textbf{Case $T_0$, (b). $\Delta^0 (x,a,b) \neq \emptyset$ and some $\theta (x,a,b)$ is between $z_s, y_{s+1}$, possibly equal to either $z_s$ or $y_{s+1}$}.

This implies that $\dist (1,z_s)\lqq \frac{1}{2} (r+p-R)$. As in Case $T_0$, (a), we may reduce the problem of estimating the operator to the case when $s$ is fixed, as well as the distances
\begin{equation}\label{eq:distan1}
 \dist(1,y_1),\dist (y_1,z_1), \dist(z_1,y_2), \dist (y_2, z_2) \ldots,\dist(z_{s-1}, y_s), \dist (y_s, z_s)\, ,
\end{equation}
 who take the respective fixed values $a_1, b_1,a_2,b_2,\dots, a_s, b_s$. As before we use the notation $A=(a_1, b_1,a_2,b_2,\dots, a_s, b_s)$ and $\Lambda (A), \Lambda_t (A)$ from \ref{notat:dmlA}.

The elements $x$ over which is taken the first sum on the left hand side of \eqref{eq:ti} can then all be decomposed as $x= h \chi \chi' k$, where $|h \chi| = \frac{1}{2}(r+p-R)$, $h\in \Lambda_t (A), \chi, \chi'\in C^{\eta}, |k| \lqq 1/2(p+R -r)$. The left-hand side of \eqref{eq:norm} is then bounded by
$$
\sum_{h\in \Lambda_t (A), \chi, \chi'\in C^{\eta}, |k| \lqq 1/2(p+R -r)}\left(\sum_{|\ell| = 1/2(r+R -p)} f(h\chi \ell ) g(\ell\iv \chi'k ) \right)^2\, ,
$$

which, by the Cauchy-Schwarz inequality, is at most
$\Pi_1 \, \Pi_2$, where
$$
\Pi_1=\sum_{h\in \Lambda_t (A), \chi \in C^{\eta}, |\ell| = 1/2(r+R -p)} [f(h \chi \ell )]^2
$$ and
$$
\Pi_2= \sum_{|\ell'| = 1/2(r+R -p), \chi'\in C^{\eta},  |k| \lqq 1/2(p+R -r)} [g(\ell' \chi' k )]^2\, .
$$

We check how many times can $[f(u)]^2$ appear in $\Pi_1$, for an arbitrary $u\in G$. If $\chi$ is non-trivial then an argument as in $T_0, (a),$ yields that it can appear only once.

If $\chi$ is trivial then the \erccs $\backslash y_s, z_s /$ of $1,x$ might be extended to an \erccs $\backslash y_s, z_s' /$ maximal, in the sense of $\Subset $, for $1,u$. The element $z_s'$ is uniquely determined for each given $u$, and the point $\theta$ must be on a geodesic $[y_s, z_s']$. By the Remark \ref{rem:stdec}, \eqref{rs2}, and Lemma \ref{lem:rt2} there are at most $4$ possible choices for $\theta$, hence for $h$.

We conclude that $\Pi_1 \lqq 4 \| f\|^2$. A similar argument shows that $\Pi_2 \lqq 4 \| g\|^2$.

\endcomment

We continue with a discussion of the cases when $\Delta^i (x,a,b) \neq \emptyset$ for $i=1,2,3,4$.

\comment
In each of these cases the central triangle intersects the side $[1,x]$:

 \begin{itemize}
   \item  either in an $\eta$--relator-tied sub-geodesic (cases $T_1, (a.I), (a,II), (a,III)$, $T_1, (b)$ and $(c)$, $T_3$ and $T_4$);
   \item or in two $\eta$--relator-tied sub-geodesics separated by a geodesic with label in $S^\eta (R)$ (case $T_2, (IV)$),
   \item or in a sub-geodesic with label in $S^\eta (R)$ (case $T_1, (a.IV)$).
 \end{itemize}

In the last case, the pair with endpoints of the intersection may be either nested in a pair $y_s, z_s$ or $z_s, y_{s+1}$, or may be separated by a point $z_s$ or $y_{s+1}$. In the other two cases, either the pair with endpoints of the whole intersection or the pair with endpoints of the first $\eta$--relator-tied sub-geodesic in the intersection is nested in a pair $y_s, z_s$ with endpoints of an \ercc. Since the number of possibilities for $s$ is linear in $r$ we may assume as before that $s$ is fixed. We may likewise assume that the distances

\begin{equation}\label{eq:distan2}
 \dist(1,y_1),\dist (y_1,z_1), \dist(z_1,y_2), \dist (y_2, z_2) \ldots,\dist(z_{s-1}, y_s)\mbox{ and if applicable }\dist (y_s, z_s)
\end{equation}
 take fixed values $a_1, b_1,a_2,b_2,\dots, a_s, b_s$ respectively. We denote
 $$A=(a_1, b_1,a_2,b_2,\dots, a_s, b_s)\, ,
 $$ or the same without $b_s$, as the case may be.

 \endcomment

\begin{cvn}\label{conv:edges}
The distances $\dist (1, C_1), \dist (a, C_a)$ are at most $r$, hence the number of possibilities for these values is linear in $r$. We therefore assume that these two lengths are fixed numbers $H, L \in [0,r]$.

All the contours that appear in the filling of the central triangle in Figures \ref{fig:bigon} and \ref{fig:triangle}, and that intersect the edge $[C_1, C_a]$, except possibly the central contour for triangles of type $T_2$ and $T_3$, and the contour appearing on the extreme right of $T_1$, have lengths at most $\frac{r}{\frac{1}{2}-2\lambda  } $. Thus, for each of the lengths of the sub-arcs of such contours, there is a number of possible choices that is linear in $r$. We therefore assume that for all such sub-arcs the lengths are fixed in $[0,r]$. In what follows for these fixed lengths we use the notation $L_i$ with $i\in \N$.
\end{cvn}

We begin with the case $T_4$, which has more similarities with the Case $T_0$.
\me


\noindent \textbf{Case $T_4$. $\Delta^4 (x,a,b) \neq \emptyset$}.

Let $x$ be a fixed element of length $p$. Every $a,b$ with $ab=x$ and $\Delta^4 (x,a,b) \neq \emptyset$ determines a decomposition of $x$ as $x= h\gamma \xi k$ with $|h|=H$, $\gamma \in RC^\eta , |\gamma | = L_1$ and $(\xi , k) \in \dd_2$ (following the notation \ref{notat:dd_i}).

\begin{figure}[htb]
\centering
\includegraphics[scale=0.8]{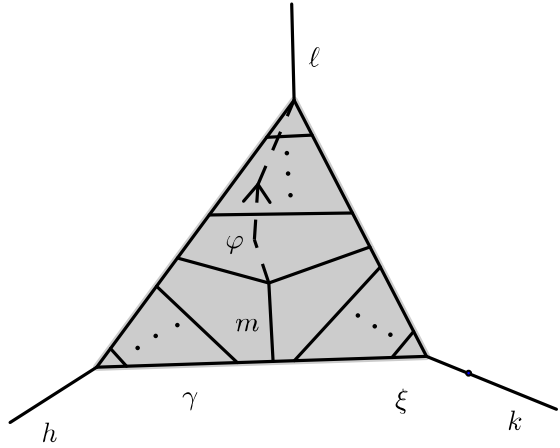}
\caption{Data describing a triple $a,b,c$ with $\Delta^4 (x,a,b) \neq \emptyset$.}
\label{fig:NotatT4}
\end{figure}

The only other data that must be added to determine the triple $x,a,b$ entirely is (see Figure \ref{fig:NotatT4} that illustrates the notation):

 \begin{itemize}
   \item  the element $m$ with $|m|=L_2$ such that $h\gamma m$ is the point in $h\G (\gamma )$ representing the median point of the central tripod;
   \item the element $\vf \in RC^\eta $ with $|\vf | = L_3$ sending that median point to $C_a$ by right translation;
   \item the element $\ell$ with $|\ell| = L$ sending $C_a$ to $a$ by right translation.
 \end{itemize}

 Thus each sum $S^4 (x)$ in \eqref{eq:norm} is taken over all possible choices of decompositions $x= h\gamma \xi k$ as described, and over all possible choices of $m, \vf$ and $\ell$.

 For a given $x$ there are two choices of $h$ and two choices of $\gamma$, and once both $h$ and $\gamma$ are chosen, $\xi$ and $k$ are uniquely defined, by the uniqueness of the \dec decomposition.

 Thus in the left hand side of \eqref{eq:norm} one may apply Cauchy-Schwartz for the first sum inside brackets and by eventually increasing the total value, remove this sum outside and add  a multiplicative factor of $4$. After that, by further increasing the value, one may replace the double sum (over all $x$ of length $p$, then for each $x$ over all decompositions $x= h\gamma \xi k$, by one sum over all $h$ of length $H$,  $\gamma \in RT^\eta , |\gamma | = L_1$ and $(\xi , k) \in \dd_2$ with $|\xi k|= p-H-L_1$ (i.e. over all products $h\gamma\xi k$). The new sum may have more terms, because some products might give elements of length $<p$.

 In other words, an upper bound of the left-hand side in \eqref{eq:norm} is
 $$
4\sum_{\substack{|h|=H,\gamma \in RT^\eta , |\gamma | = L_1,\\[2pt] (\xi , k) \in \dd_2 , |\xi k|= p-H-L_1} }\left( \sum_{\vf \in RC^\eta,|\vf | = L_3, |\ell| = L } f(h\gamma m\vf\ell ) g(\ell^{-1} \vf\iv m\iv \xi k) \right)^2\!.
 $$
 The Cauchy-Schwarz inequality applied to the sum in between round brackets yields the upper bound
 \begin{equation}\label{eq:CS1}
 4 \sum_{h, \gamma, (\xi , k) }\left(\sum_{\vf \in RC^\eta,|\vf | = L_3, |\ell| = L } [f(h\gamma m\vf\ell )]^2 \right) \left( \sum_{|\ell'| = L, \vf'\in RC^\eta,|\vf'| = L_3 } [g(\ell' \vf' m\iv \xi k)]^2 \right)\!.
 \end{equation}

 In the sums above, the element $m$ is completely determined by the choice of the two elements $\gamma $ and $\xi \, $.

 Let $M(\gamma )$ be the set of two elements $m$ in $G$ such that $\gamma m$ is a point in $\G (\gamma )$ at distance $L_2$ from $\gamma$. Likewise let  $M(\xi )$ be the set of two elements $m\iv$ in $G$ such that $m$ is a point in $\G (\xi )$ at distance $L_2$ from the identity element $1$.

 We increase the two sums that appear in brackets in \eqref{eq:CS1} by adding a summation over $m\in M(\gamma )$ for the first, and a summation over $m\in M(\xi )$ for the second. We obtain an upper bound of the form $4\Pi_1 \Pi_2$, where
 $$
 \Pi_1 = \sum_{|h|=H, \gamma \in RT^\eta , |\gamma | = L_1,m\in M(\gamma ), \vf\in RC^\eta,|\vf| = L_3, |\ell| = L } [f(h\gamma m \vf \ell )]^2
 $$ and
$$
 \Pi_2 = \sum_{|\ell'| = L, \vf'\in RC^\eta,|\vf' | = L_3, (\xi , k) \in \dd_2, |\xi k|= p-H-L_1, m\in M(\xi) } [g(\ell' \vf' m \xi k)]^2\, .
 $$

For each $u\in G$ we check how many times can $[f(u)]^2$ appear in $\Pi_1$. Corollary \ref{cor:intermediate} implies that for the given $u$ there are at most two possible choices for each of the factors $h$, $\gamma ,\, \ell$ and $\vf$. This implies that each  $[f(u)]^2$ can appear in $\Pi_1$ at most $2^4$ many times, whence $\Pi_1 \lqq 2^4 \| f\|^2\, $.

For an arbitrary $v\in G$ we count the number of appearances of $[g(v)]^2$ in $\Pi_2$. By  Corollary \ref{cor:intermediate} there at most two possible choices for $\ell'$, respectively for $\vf'$. Once $\vf'$ fixed, $m$ is such that $\ell' \vf'm$ is a point on the contour that comes next to the one containing $\ell' \vf'$ in $\G (v)$, and at distance $L_2$ from  $\ell' \vf'\, $. Therefore there are at most $2$ possible choices for $m$.

Once $\ell', \vf'$ and $m$ are fixed, $(\xi , k)$ are uniquely determined, since $\xi$ is the first element in the \dec decomposition of $(\ell' \vf' m)\iv u$, as defined in Corollary \ref{cor:intermediate}.

We may then conclude that $\Pi_2 \lqq 2^3 \| g\|^2\, $.

Thus we obtain that
$$
\Pi_1 \Pi_2 \lqq 2^7 \|f\|^2 \|g\|^2\, .
$$

\medskip

\noindent \textbf{Case $T_2$. $\Delta^2 (x,a,b) \neq \emptyset$}.

We denote by $\omega$ the central contour of the central triangle, and by $\omega_{\alpha\beta}$ the intersection $\omega \cap [C_\alpha , C_\beta ]$. We denote by $\vartheta, \rho ,\zeta $ the element in $G$ that is defined respectively by the labels of the paths $\omega_{1x}, \omega_{1a}, \omega_{ax}$. We denote by $\epsilon_\alpha$ the label of the intersection of $\omega$ with the corner of $\alpha \in \{1,a,x\}$.

\begin{figure}[htb]
\centering
\includegraphics[scale=0.8]{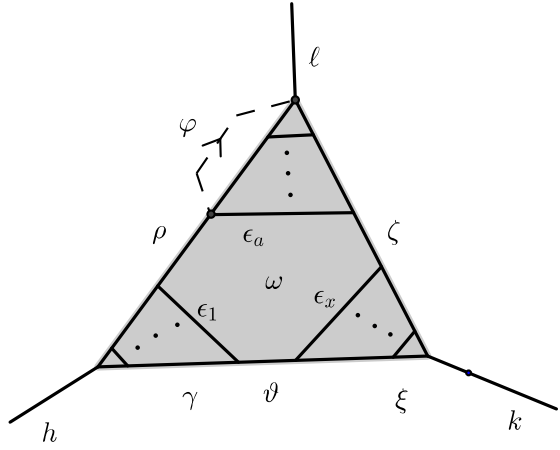}
\caption{Data describing a triple $a,b,c$ with $\Delta^2 (x,a,b) \neq \emptyset$.}
\label{fig:NotatT2}
\end{figure}

We denote by $\vf $ the label of the geodesic joining the upper corner of $\omega_{1a}$ with $C_a$. The rest of the notations are identical to the ones in Case $T_4$. See Figure \ref{fig:NotatT2}.

As in Case $T_4,$ all the contours filling the central triangle and intersecting $[C_1, C_a]$ have length at most $\frac{r}{\frac{1}{2} - 2 \lambda }$, except possibly $\omega$.

Therefore, without loss of generality, we assume that all the lengths of sub-arcs of such contours appearing as intersections with other contours or with the edges of the geodesic triangle are fixed. Besides the lengths for which notations were already introduced in Case $T_4$, the other fixed lengths are denoted by $L_i$ with $i\geqslant0$.

For each $x$ of length $p$, the sum $S^2 (x)$ is first over the possible decompositions $x= h \gamma \vartheta \xi k$, then over all possible $\omega , \vf$ and $\ell\, $.

For a given $x$, when counting the number of decompositions $x= h\gamma \vartheta \xi k\, ,$ one notes that there are $2^2$ possible choices for the pair $(h, \gamma )$ since the lengths of both components are fixed. Once these two are chosen, $\vartheta$ labels a sub-geodesic in $\mathcal{G}(x)$ entirely contained in one contour $t_i$ and joining the point $h\gamma$ to the nearest endpoint of the intersection $t_i \cap t_{i+1}$, where $t_{i+1}$ is the contour in  $\mathcal{G}(x)$ consecutive to $t_i$. Thus, once $h$ and $\gamma$ are chosen, the choice of $\vartheta $ is compulsory. And once $h, \gamma$ and $\vartheta$ are given, $\xi$ and $k$ are uniquely determined. Thus we may obtain an upper bound with the sum over the decompositions removed outside the round bracket of \eqref{eq:norm}, and a multiplicative factor of $2^2$ added. By further increasing the bound, we may consider from the beginning the sum over all possible $h, \gamma, \vartheta , \xi $ and $ k$.

We thus obtain that the left hand side of \eqref{eq:norm} is bounded by
\begin{equation}\label{eq:boundT2}
2^2 \sum_{h, \gamma , \vartheta, \xi , k} \left( \sum_{\omega , \vf \in RC^\eta, |\vf |=L_4, |\ell | =L } f(h\gamma \epsilon_1 \rho \vf \ell ) g(\ell\iv \vf \iv \epsilon_a \zeta \epsilon_x \xi k) \right)^2\, .
\end{equation}

\medskip

\noindent \textbf{Case $T_2$, (I).} $|\omega_{1x}| < \lambda |\omega|$. This implies that $|\omega_{1a}|, |\omega_{ax}| > \left(\frac{1}{2}-4\lambda \right) |\omega|$.

In particular, it follows that $|\omega | \leqslant\frac{r}{\frac{1}{2}-4\lambda}$, which implies that we may assume without loss of generality that all the $\omega_{\alpha\beta}$ and $\epsilon_{\alpha}$ have fixed lengths.

Each decomposition $ab=x$ having a central triangle with the properties described above determines a decomposition of $x$ as $x= h \gamma \vartheta \xi k$, where $|h|=H, \gamma \in RC^\eta , |\gamma |=L_1, \vartheta \in S (R), |\vartheta |=L_2 $ and $(\xi , k ) \in \dd_2$, with the notation \ref{notat:dd_i}, $\xi k$ of length $p-H-L_1-L_2$.

In order to completely describe the decomposition $ab=x$ the data that must be added to the above is:
\begin{itemize}
  \item the contour $\omega$ with a sub-geodesic labeled by $\vartheta$;
  \item the element $\vf \in RC^\eta $ sending the upper corner of $\omega_{1a}$ to $C_a$ by right translation; this element has fixed length $L_4$;
  \item the element $\ell$ of fixed length $L$ sending $C_a$ to $a$ by right translation.
\end{itemize}

We apply the Cauchy-Schwarz inequality to \eqref{eq:boundT2} and obtain the upper bound
$$
2^2 \sum_{h, \gamma , \vartheta, \xi , k} \left[ \sum_{\substack{\omega \supseteq \vartheta, \vf \in RC^\eta,\\ |\vf |=L_4, |\ell | =L} } f(h\gamma \epsilon_1 \rho \vf \ell )^2 \right] \left[ \sum_{\substack{\omega \supseteq \vartheta, \vf' \in RC^\eta,\\ |\vf' |=L_4, |\ell' | =L }}g(\ell' \vf'  \epsilon_a \zeta \epsilon_x \xi k)^2 \right].
$$

In the above the notation $\omega \supseteq \vartheta $ signifies that the contour $\omega$ has the sub-path following its intersection with $\mathcal{G}^\eta (\gamma )$ and of the fixed required length labeled by the element $\vartheta$.

This gives the upper bound $2^2\Pi_1\Pi_2$, where
$$
\Pi_1= \sum_{h, \gamma \in S^\eta (R), \omega , \vf \in RC^\eta, |\vf |=L_4, |\ell | =L } [f(h\gamma \epsilon_1 \rho \vf \ell )]^2
$$ and
$$
\Pi_2= \sum_{|\ell' | =L, \vf' \in RC^\eta, |\vf'|=L_4, \omega , \xi \in RC^{\eta}} [g(\ell' \vf' \epsilon_a \zeta \epsilon_x \xi k)]^2\, .
$$

At this point in the argument, we wish to emphasize that the condition $\lambda \lqq \frac{1}{10}$ is necessary to be able to proceed. Indeed, if in each of the sums above we may replace the summation over all the possible choices for $\omega$ by respective summations over all the possible choices for $\epsilon_1$, $\rho$, respectively $\epsilon_a, \zeta , \epsilon_x  $, all of fixed lengths, then the desired inequality can be easily obtained, as shown in the sequel. Both replacements can be done if we know that both $|\omega_{1a}|$ and $ |\omega_{ax}|$ are at least $ \lambda |\omega|$. Since the only hypothesis at our disposal is that $ |\omega_{1a}|, |\omega_{ax}| > \left(\frac{1}{2}-4\lambda \right) |\omega|$, the above can be granted only by $\lambda \lqq \frac{1}{10}$.

{If, on the other hand, for at least one of the two lengths $|\omega_{1a}|$ and $ |\omega_{ax}|$ we have that it is $<\lambda |\omega|$, say $ |\omega_{ax}|$, then in the second sum the same term $[g(\ell' \vf' \epsilon_a \zeta \epsilon_x \xi k)]^2$ appears at least as many times as the number of contours $\omega$ such that $\epsilon_a \zeta \epsilon_x \sqsubset \mathrm{label} (\omega )$. In general, there may be infinitely many such $\omega$ for one fixed term, and even if we assume moreover that $|\omega_{1a}| \gq \lambda |\omega|$, which implies that $|\omega| \lqq \frac{r}{\lambda}$, we can still have $\exp (\alpha r)$ distinct contours $\omega$ such that $\epsilon_a \zeta \epsilon_x \sqsubset \mathrm{label} (\omega )$, for some $\alpha >0$.

Thus, we do need that both $|\omega_{1a}|$ and $ |\omega_{ax}|$ are at least $ \lambda |\omega|$, hence that $\lambda \lqq \frac{1}{10}$. This allows to replace in $\Pi_1$ the summation over $\omega $ by a summation over $\rho$, and in $\Pi_2$ the summation over $\omega $ by a summation over $\zeta$.

The factor $\Pi_1$ is then bounded by $2^5 \| f\|^2$ because the lengths of all factors $h,\gamma , \epsilon_1,$ $\rho, \vf ,\ell$ are fixed, hence every element $y\in G$ can be written as a product $h\gamma \epsilon_1 \rho \vf \ell$ in at most $2^5$ different manners, by Corollary \ref{cor:intermediate}.

Now consider $\Pi_2$. Here too, since the lengths of $\ell' ,\vf' ,\epsilon_a, \zeta , \epsilon_x$ are fixed, for a given $y\in G$ there are only at most $2^5$ different choices possible for these elements so that their product is a prefix of $y$. Once these choices are made, the pair $(\xi ,k )$ is uniquely determined.

\medskip

\noindent \textbf{Case $T_2$, (II)} $|\omega_{1x}| \gq \lambda |\omega|$ and $|\omega_{1a}| \gq \lambda |\omega |$.

In this case, the choice of $\vartheta $ entirely determines the choice of $\omega$, of the elements  $\epsilon_\alpha$ and $\rho, \zeta$; likewise, the choice of $\rho $ entirely determines the choice of $\omega$, of the elements  $\epsilon_\alpha$ and $\vartheta , \zeta$. Also, since $|\omega | \lqq \frac r\lambda$ it follows again that one can reduce to the case when all the intersection sub-arcs contained in $\omega$ have fixed lengths.

We apply the Cauchy-Schwarz inequality again for the sum inside the round brackets in \eqref{eq:boundT2}. The sum in $g$ will depend on $\ell, \vf , \xi , k$ and to make it independent of $\vartheta $, we add a further sum over all possible $\epsilon_a, \zeta$ and $\epsilon_x$ of that fixed length. After that, the upper bound can be separated into a product $2^2\Pi_1\Pi_2$, where
\begin{equation}\label{eq:pi1}
\Pi_1= \sum_{|h|=L_1, \gamma \in RT^{\eta }, |\gamma |=L_2,  \rho \in S^{\lambda } (R), |\rho |= L_3, \vf \in RC^\eta, |\vf |=L_4, |\ell | =L_5} [f(h\gamma \epsilon_1 \rho \vf\ell )]^2
\end{equation}
 and
 \begin{equation}\label{eq:pi2}
  \Pi_2= \sum_{|\ell' | =L_5, \vf'\in RC^\eta, |\vf '|=L_4, \zeta\in S (R), |\zeta |= L_6, |\epsilon_a|=L_7, |\epsilon_x| = L_8 ,  (\xi , k)\in \dd_2 } [g(\ell' \vf ' \epsilon_a \zeta \epsilon_x \xi k)]^2\, .
 \end{equation}

For a given $u\in G$ when one looks for the number of decompositions $u= h\gamma \epsilon_1 \rho \vf \ell$, and since all the lengths appearing are fixed, a bound by $2^5\| f\|^2$ is obtained as previously.

Likewise, for the second factor one studies for an arbitrary element $v\in G$ the number of decompositions $v= \ell' \vf' \epsilon_a \zeta \epsilon_x \xi k$. The lengths of $\ell', \vf', \epsilon_a, \zeta, \epsilon_x $ being fixed, the number of distinct choices for a fixed $v$ is at most $2^5$. Once all are fixed, $\xi$ is uniquely determined as first component in the decomposition of $\xi k$ defined in Theorem \ref{cor:sdecomp}.

\medskip

\noindent \textbf{Case $T_2$, (III).} $|\omega_{1x}| \gq \lambda |\omega|$ and $|\omega_{1a}|< \lambda |\omega |$. The latter inequality implies that in fact both $|\omega_{1x}| $ and $ |\omega_{ax}| $ are at least $\left( \frac{1}{2} - 4 \lambda \right) |\omega |$, in particular $|\omega_{ax}| \gq \lambda |\omega|$ as $\lambda \leqslant 1/10$.

In this case, the respective lengths of $\omega , \vartheta , \zeta $ and $\epsilon_x$ cannot be fixed in terms of $r$. The length of $\rho$ can be assume fixed, because it is at most $r$ hence the number of possibilities for it is linear in $r$. The lengths of $\epsilon_1$ and $\epsilon_a$ can likewise be fixed, because they are intersections of $\omega$ with cells of lengths $O(r)$.

Moreover, $\omega $ and $ \vartheta$ entirely determine each other, as well as the $\epsilon_\alpha$ with $\alpha \in \{1,a,x\}\, $, and $\zeta, \rho$.

We apply the Cauchy-Schwarz inequality to the sum inside the round brackets in  \eqref{eq:boundT2}. In order to make the sum in $f$ independent of $\vartheta$, we increase it, by adding an extra sum over all $\rho \in S(R)$ of fixed length $L_3$ and all $\epsilon_1\in S(R)$ of length $L_5\, $.

We obtain the upper bound $2^2\Pi_1 \Pi_2$ with
\begin{equation}\label{eq:pi1a}
\Pi_1= \sum_{|h|=L_1, \gamma \in RT^{\eta }, |\gamma |=L_2, |\epsilon_1|=L_5,  \rho \in S (R), |\rho |= L_3, \vf \in RC^\eta, |\vf |=L_4, |\ell | =L} [f(h\gamma \epsilon_1 \rho \vf\ell )]^2
\end{equation}
 and
 \begin{equation}\label{eq:pi2a}
  \Pi_2= \sum_{|\ell' | =L, \vf'\in RC^\eta, |\vf '|=L_4, |\epsilon_a|=L_6, \zeta\in S^\eta (R), \epsilon_x ,  (\xi , k)\in \dd_2 } [g(\ell' \vf ' \epsilon_a \zeta \epsilon_x \xi k)]^2\, .
 \end{equation}

For the estimate of $\Pi_1$ this makes no difference, since the only important thing was that the lengths of the elements intervening in the product $h\gamma \epsilon_1 \rho \vf\ell$ were fixed.

For an arbitrary $v$, in a decomposition $v= \ell' \vf ' \epsilon_a \zeta \epsilon_x \xi k$, there are 8 possible choices for $\ell', \vf' , \epsilon_a$. Once these are fixed there is a unique choice for the contour whose label would contain $\zeta$, hence two possible choices for $\zeta , \epsilon_x$. Once these are likewise chosen, there is only one possibility for the pair $(\xi , k)$.

\medskip

\noindent \textbf{Case $T_3$. $\Delta^3 (x,a,b) \neq \emptyset$}. This may be seen as a particular sub-case of Case $T_2$, in which either $\vartheta, \zeta$ or $\rho$ are trivial. The respective sub-cases that cover each of these situations are $T_2, (I)\, ,$ $T_2, (II)\, $ and  $T_2, (III)\, $. In each of these cases, the possibility that $\vartheta, \zeta$ or $\rho$ equal $1$ is not excluded.

\medskip

\noindent \textbf{Case $T_1$, (a). $\Delta^1 (x,a,b) \neq \emptyset$ and the side contained in one contour is $[C_1 , C_x]$}.

This can be seen as a particular sub-case of $T_2$, with $\gamma = \xi = \epsilon_1 = \epsilon_x =1\, $. Therefore we keep a consistent notation: we denote by $\omega$ the last contour in the central triangle, and by $\omega_{\alpha \beta}$ its intersection with $[C_\alpha , C_\beta]$, where $\alpha, \beta \in \{1,a,x \}$. We denote the elements in $G$ represented by the labels of $\omega_{1x}\, ,\, \omega_{1a}$ and $\omega_{ax}$ by $\vartheta, \rho $ and $\zeta$, respectively. We denote by $\epsilon_a$ the element of $G$ represented by the label of the intersection of $\omega$ with the contour stacked above him.

Let $\vf$ be the element right-translating the upper endpoint of $\omega_{1a}$ to $C_a$; let $\ell$ be the element right-translating $C_a$ to $a$. We record below the differences in the argument in this case, compared to Case $T_2$.

As previously, we may assume that the lengths $H$ of $h$, $L$ of $\ell$, as well as the lengths of $\rho$,  $\vf $, $\epsilon_a$, are fixed and $O(r)$.

Given an element $x$ of length $p$, the sum in $S^1 (x)$ can be written as a double sum, first over all decompositions $x= h\vartheta k$, then over all possible $\omega$ (if $\omega$ is not already completely determined by $\vartheta$), $\vf$ and $\ell$.

Assume that $|\omega_{1x}| \gq \lambda |\omega |$. In this case for a given $x$ there are two possible choices for $h$ (of fixed length), and once $h$ is chosen, there is only one possible choice for the contour containing the sub-geodesic labeled by $\vartheta$, hence there are only two possible choices for $\vartheta$. In total each $x$ has at most $2^2$ decompositions in this case.

Assume that $|\omega_{1x}| < \lambda |\omega |$. This implies that $\omega_{1a}, \omega_{ax}$ have lengths at least $\left(\frac{1}{2}- 2\lambda \right)|\omega|$, in particular $\omega $ has length at most $\frac{r}{\frac{1}{2}- 2\lambda}$. Consequently, we may argue that without loss of generality the lengths  of $\vartheta$ and $\zeta$ may be assumed fixed as well. Therefore, given an element $x$, the number of possible decompositions $x= h\vartheta k$ is at most $2^2$.

We therefore argue as in the previous cases and find that an upper bound for the left hand side in \eqref{eq:norm} is
$$
4 \sum_{|h|=H, \vartheta\in S(R) , k}\left(\sum_{\omega , |\vf|=L_4\, ,\, |\ell|=L} f(h\rho \vf \ell ) g(\ell\iv \vf\iv \epsilon_a \zeta k) \right)^2\, .
$$

Case $T_1$, (a), (I), $|\omega_{1x}| <\lambda |\omega |$, is discussed exactly like Case $T_2$, (I). Likewise, Case $T_1$, (a), (II), $|\omega_{1x}|\gq \lambda |\omega |,\, |\omega_{1a}|\gq \lambda |\omega |$, is treated like Case $T_2$, (II); Case $T_1$, (a), (III), $|\omega_{1x}|\gq \lambda |\omega |,\, |\omega_{1a}|<\lambda |\omega |$, is treated like Case $T_2$, (III).

\medskip

\noindent \textbf{Case $T_1$, (b). $\Delta^1 (x,a,b) \neq \emptyset$ and the side contained in one contour is $[C_1 , C_a]$}.

This can be identified with a particular sub-case of Case $T_2$, in which $\gamma = \epsilon_1=\epsilon_a = \vf =1$. We therefore keep the same notations and conventions as in Case $T_2$.

For the very first step of the argument, in order to obtain an upper bound of type \eqref{eq:boundT2} we must see how may different decompositions of the form $x= h \vartheta \xi k$ can a fixed element $x$ have. Since $h$ has fixed length, there are two possible choices for it.

If $\vartheta$ has length at least $\lambda |\omega |$ then once $h$ chosen, there is a unique possibility for $\omega$ and two possible choices for $\vartheta$.

If $\vartheta$ has length $<\lambda |\omega |$, then both $\rho$ and $\zeta$ have length at least $\left( \frac{1}{2}-2\lambda \right)|\omega|$. In particular $|\omega |\lqq \frac{r}{\frac{1}{2}-2\lambda}\, $, and without loss of generality we may assume that $\vartheta , \rho , \zeta , \epsilon_x$ have fixed length.

In this case also there are two possible choices for $\vartheta$. Once $h$ and $\vartheta$ chosen, $(\xi, k)$ is uniquely determined.
Thus, we again obtain a bound as in \eqref{eq:boundT2}.

In Case $T_1$, (b), (I), $|\omega_{1x}| <\lambda |\omega |$, the argument follows the one in Case $T_2$, (I). The only difference is that $\Pi_1$ has as an upper bound $2^2 \|f\|^2\, ,$ and $\Pi_2$ is at most $2^2 \|g\|^2\, .$

In Case $T_1$, (b), (II), when $|\omega_{1x}| \gq \lambda |\omega |$ and  $|\omega_{1a}| \gq \lambda |\omega |$, the discussion is again similar to the one in Case $T_2$, (II); in the end we obtain that $\Pi_1\lqq 2^2 \|f\|^2\, ,$ and $\Pi_2\lqq 2^3 \|g\|^2\, .$

The discussion in Case $T_1$, (b), (III), when $|\omega_{1x}| \gq \lambda |\omega |$ and  $|\omega_{1a}| < \lambda |\omega |$, is likewise similar to the one in Case $T_2$, (III), with slight modifications of the multiplicative factors in the end.

\medskip

\noindent \textbf{Case  $T_1$, (c). $\Delta^1 (x,a,b) \neq \emptyset$ and the side contained in one contour is $[C_a , C_x]$}.

This again may be seen as a sub-case of $T_2$, in which $\xi = \vf = \epsilon_a = \epsilon_x =1$. We keep the rest of the notation accordingly.

\me
 To obtain the first estimate of the left hand side of \eqref{eq:norm} by a sum as in \eqref{eq:boundT2}, we argue very similarly to Case  $T_1$, (b). Then we consider three distinct sub-cases of $T_1$, (c), denoted (I), (II), and (III), defined by the same inequalities as cases $T_2$, (I), (II), and (III) respectively, and discussed in the same way, with the appropriate changes in the multiplicative factors.  \endproof

\section{ Applications of the RD property for $C'(1/10)$--groups.}\label{sec:mainappl}

\subsection{Approximation properties}\label{sec:applapprox}

We have proven that every finitely generated $C'(1/10)$--group has property RD with respect to a word length function.
We now combine this result with an older theorem in order to get Corollary~\ref{cor:map}, that is the Metric Approximation Property for the reduced $C^*$--algebra and the Fourier algebra of such a group.
We refer to \cite{Eymard, JolissaintValette:RDapp} for a definition of the Fourier algebra.

\begin{thm}(\cite{Haagerup:RD_F, JolissaintValette:RDapp,BrodNiblo:approximation})\label{thm:approx}
If a discrete group $G$ is endowed with a length function defining a pseudo-metric that is conditionally negative definite, and it has the RD property with respect to this length function then the reduced $C^*$--algebra $C^*_r(G)$ and the Fourier algebra $A(G)$ have the  Metric Approximation Property.
\end{thm}

 In general, a word length metric is not conditionally negative definite. Fortunately, for the groups that we consider there exists another metric
 which is conditionally negative definite and bi-Lipschitz equivalent to every word metric.

  Haglund and Paulin defined  the concept of space with walls \cite{HaglundPaulin}. Such a space is naturally endowed with the so-called wall pseudo-metric, which is conditionally negative definite (see for instance \cite[Theorem 12.2.9]{BrownOzawa:book}). It therefore suffices to find on every finitely generated
$C'(1/10)$--group $G$ a structure of space with walls with pseudo-metric bi-Lipschitz equivalent to a word metric. In what follows we briefly explain how the construction of such walls, done by Wise in \cite{Wise:cancell} for small cancellation groups with finite presentation, extends to infinite presentations.  For the bi-Lipschitz equivalence of the wall metric with a word metric we use the \dec decomposition, Theorem~\ref{lem:intermediate}. Note that Wise's construction can be performed under small cancellation conditions weaker than $C'(1/8)$, still for our purposes arguing under the condition $C'(1/8)$ suffices.

We refer the reader to~\cite{ArzhOsajda:Haagerup} for a stronger result, which shows that a large class of small cancellation groups, including infinitely presented $C'(1/6)$-groups,  have a discrete structure of walls whose pseudo-metric is bi-Lipschitz equivalent to a word metric.

\me

\noindent (I) Let $X$ be the $2$--dimensional complex obtained from the Cayley graph of $G$ by glueing a $2$--cell along any cycle with boundary labeled by a relator. Without loss of generality one may assume that all $2$--cells have as boundaries polygons with an even number of edges, otherwise one may consider the barycentric subdivision for every $1$--cell.

\me

\noindent (II) On each $2$--dimensional cell, the mid-points of opposite edges are joined by a new edge, which we call for simplicity \emph{median edge}. A maximal connected union of median edges intersecting each $2$--cell in at most one such edge is called a \emph{hypergraph}.
For each hypergraph, the union of $2$--cells intersecting it in a median edge is called a \emph{hypercarrier}.

\me

\noindent (III) If the presentation of $G$ satisfies $C'(1/6)$, or even weaker small cancellation properties \cite{Wise:cancell}, then each hypergraph is a tree embedded into the complex $X$ \cite[Corollary3.12]{Wise:cancell}, which separates $X$ into exactly two connected components, one for each of the two connected components of the boundary of the hypercarrier \cite[Lemma 3.13]{Wise:cancell}. Note that in \cite{Wise:cancell} the group $G$ is supposed to have a finite presentation. This is not necessary for this part of the argument: given a group $G$ presented as in \eqref{eq:pres} and satisfying the $C'(1/6)$--condition, $G$ is a direct limit of the finitely presented groups $G_k$ described in \eqref{eq:presK}. The above statements are true for each 2-complex $X_k$ corresponding to each $G_k.$ Therefore, they are true for the complex $X$ corresponding to $G$: the hypergraphs in $X$ are increasing countable unions of simplicial trees, therefore they are simplicial trees, and the separation properties easily follow from those in $X_k$ since every path joining two points in $X$ is covered by finitely many $2$-cells, hence a copy of it already appears in a complex~$X_k$.

\me

\noindent (IV) Let $\eta \gq \frac{3}{8}$ and let $a$ and $b$ be two arbitrary elements in $G$. According to Theorem \ref{lem:intermediate}, the set $\mathcal{G}^\eta (a,b)$ contains every geodesic with endpoints $a,b$. It follows that every hypergraph separating $a$ and $b$ must intersect $\mathcal{G}^\eta (a,b).$ That is,  it must contain at least one midpoint of one edge in $\mathcal{G}^\eta (a,b)$. Conversely, every hypergraph containing a midpoint of one edge in $\mathcal{G}^\eta (a,b)$ separates $a,b$, since clearly $a$ and $b$ are each connected to a different connected component of the boundary of the hypercarrier of that hypergraph.

      It remains to note that the number of edges in $\mathcal{G}^\eta (a,b)$ is bounded by a multiple of $\dist (a,b)$. With the notation from Theorem \ref{lem:intermediate}, we have
      $$\dist (a, y_1)+ \sum_{i=2}^m \dist (z_{i-1}, y_i) + \dist (z_m, b) \leqslant \dist (a,b)\, .
      $$

      Also, each contour $t^{(i)}_j$ intersects a geodesic joining $a,b$ in a segment of length at least $\left( \eta - \frac{1}{4} \right)\left| t^{(i)}_j \right|\gq \frac{1}{8} \left| t^{(i)}_j \right|$ with interior disjoint of all intersections with other contours. Therefore, the sum of all the lengths of all the $t^{(i)}_j$ is at most $8 \dist (a,b)$.

This shows that the number of hyper-graphs separating the two points $a,b$ (which equals the wall pseudo-distance from $a$ to $b$) is at least $\dist (a,b)$ and at most $4 \dist (a,b)\, $.
Thus, the wall pseudo-metric is bi-Lipschitz equivalent to the word length metric, as required.

\subsection{Operator growth series}\label{sec:applgrowth}

Given a finitely generated group $G$ and a fixed finite set of generators $A$ for it, one can study the growth of $G$ with respect to $A$ by considering the \emph{spherical growth function} $\mathfrak a  :\N \to \N$ defined so that ${\mathfrak a}(n)$ is the set of all elements $g\in G$ with $|g|_A= n\, $.
A relevant object  is the \emph{spherical growth series} defined by
$$
f(z)= \sum_{n=0}^\infty {\mathfrak a} (n) z^n\, .
$$ It has been thoroughly studied for important classes of groups such as Coxeter groups \cite{Paris:growthCoxeter, Wagreich:growth,ScottR:coxeter}, Kleinian and Fuchsian groups \cite{CannonWagreich, Wagreich:growth,FloydPlotnick}, nilpotent groups \cite{Benson:abelian,Benson:growthnil,shapiro:growth}, and hyperbolic groups \cite{Cannon:graph,Gromov(1987)}.

The \emph{complete growth series} of the group $G$ is a power series with coefficients in the group ring $\Z [G ]$ defined by
$$
F(z)= \sum_{n=0}^\infty A_n z^n\, ,
$$ where $A_n$ is the sum of all the elements of length $n$ in $G$. It was defined and studied for hyperbolic groups \cite{GrigorchukNagnibeda} and for direct, free and graph products \cite{ACFMR}.

The complete growth series has been further generalized as follows \cite{GrigorchukNagnibeda}: given a representation $T: G \to \mathcal{B}_\ast$ of the group $G$ in the group of invertible elements of a Banach algebra, the \emph{operator growth series with respect to }$T$ is the power series with coefficients in $\mathcal{B}_\ast$ defined by
$$
F_T(z) = \sum_{g\in G} T(g) z^{|g|} = \sum_{n=0}^\infty A_n^T z^n \, .
$$

In the particular case of the representation of $G$ inside the reduced $C^*$--algebra of $G$, the series above is simply called \emph{operator growth series}.
For more details on these growth series we refer the reader  to \cite[Chapter VI]{delaHarpe:topics}.

Grigorchuk and Nagnibeda carried out a comparison between the series above \cite{GrigorchukNagnibeda}. They proved that:

\begin{enumerate}
  \item\label{gr1} The radius of convergence of the complete growth series equals that of the spherical growth series.

 \me

  \item\label{gr2} If a group is amenable then the radius of convergence of the operator growth series equals that of the spherical growth series.
\end{enumerate}

They conjectured that the converse of \eqref{gr2} is also true. In this context, they proved that if a group has property RD then the radius of convergence of the operator growth series equals the square root of the radius of convergence of the spherical series. In particular, this is true for the groups that we consider, see Corollary \ref{cor:opg}.

\subsection{Ergodic properties in quantum dynamical systems}\label{sec:quantum}

The property of Rapid Decay is also useful for generalizations of classical ergodic results to the setting of quantum dynamical systems (also called $C^*$--dynamical systems). This holds, in particular, for the equivalent for quantum dynamical systems of unique ergodicity: given an automorphism $\alpha$ of a unital $C^*$--algebra $A$, for every $a\in A$ the sequence of ergodic
averages $\frac{1}{n}\sum_{k=1}^n \alpha^k (a)$ converges in norm to a scalar multiple $E(a)$ of $1$ (the conditional expectation).

It was proved by Abadie and Dykema \cite{AbadieDykema} that if $G$ is a group with the property RD, then the unique ergodicity holds for the action on $C^*_{\rm{red}} (G)$ of an automorphism $\beta$ induced by an automorphism of the group $G$ with certain properties.

A stronger property of strict weak mixing for the same quantum dynamical systems was investigated in \cite{FidaleoMukh}. Moreover, in \cite{Fidaleo} an even stronger ergodic property is studied: the convergence to the equilibrium, implying all the ergodic properties mentioned above. This latter property is specific to quantum dynamical systems, and has no counterpart in the classical case. It is proven in\cite{Fidaleo} that if a group $G$ has property RD  then the action of an automorphism $\beta$ as above on the algebra $C^*_{\rm{red}} (G)$ of $G$  has the property of convergence to the equilibrium.

\section{Quasi-homomorphisms on small cancelation groups}\label{sec:qm}

Recall that a \emph{quasi-homomorphism}  (also called a \emph{quasi-morphism} and a \emph{pseudo-character}) on a group $G$  is a function $\fh :G \to \R$ such that its \emph{defect}
$$
\dgot (\fh ):= \sup_{a,b\in G}\left| \fh (ab) - \fh (a) - \fh (b) \right|
$$
is finite. The real vector space $\mathcal{Q} (G)$ of all quasi-homomorphisms of $G$ has three important subspaces: the subspace $\ell^{\infty}(G)$ of bounded real functions on $G$,
the subspace $\rm{Hom} (G,\R)= H^1 (G,\R )$ of homomorphisms on $G$, and the subspace $\ell^{\infty}(G) + \rm{Hom} (G, \R)$ of the functions that differ from a homomorphism by a bounded function.
Consider the quotient spaces
$$
QH(G) = \mathcal{Q} (G)/\mathcal{B}(G)\; \mbox{ and\,\, } \widetilde{QH}(G) = \mathcal{Q} (G)/\left[\ell^{\infty}(G) + \rm{Hom} (G, \R)\right]\,.
$$
The space $\widetilde{QH}(G)$ can be identified with the kernel of the comparison map
$$
H^2_b(G) \to H^2 (G)\, ,
$$ where $H^2_b(G)$ is the second bounded cohomology of $G$.

In this paper, as a second application of our results on the geometry of small cancelation groups with $C' \left(1/12 \right)$--condition, we show that for such a group $G$ the space $\widetilde{QH}(G)$ is infinite dimensional, with a basis of power continuum.

Following the work of Epstein and Fujiwara \cite{Epstein-Fujiwara, Fujiwara:BC1, Fujiwara:BC2} as well as of Bestvina and Fujiwara \cite{BestvinaFujiwara}, we shall prove the following.

 \begin{prop}\label{prop:hn}
 Let $G$ be a finitely generated infinitely presented group and let $\la S \mid R \ra$ be a presentation such that $R$ satisfies  $C' \left(1/12 \right)$--condition.
 For a given $\eta \in [3\lambda, \frac{1}{2}-2\lambda ]$ appropriately chosen, there exists a sequence $\fu_n$ of elements in $G$ and a sequence $\fh_{\fu_n} :G\to \R $ of quasi-morphisms, with $n\in \N\, ,\, n\geqslant1\, ,$ such that
\begin{enumerate}
  \item\label{h1} the set of word lengths $|\fu_n|$ diverges to $\infty$;
  \item\label{h2} every group homomorphism $\phi :G \to \R $ has the property that $\phi (\fu_n)=0$ for every $n\in \N\, ,\, n\geqslant1$;
  \item\label{h3} the sequence of defects $\dgot \left (\fh_{\fu_n} \right)$ is bounded;
  \item\label{h4} for every $n$ and every $k\in \N, k\geqslant1$, $\fh_{\fu_n}\left( \fu_n^k \right)=k\, $;
  \item\label{h5} for every $n\neq m$, and every $k\in \N, k\geqslant1$, $\fh_{\fu_n}\left( \fu_m^k \right)=0\, $.
\end{enumerate}
\end{prop}

\proof  We enumerate the relators $\{ r_1,r_2,\dots \} $
 in $R$ so that their lengths compose a non-decreasing sequence.
Consider the sequence of finite subsets of $\N$ defined by
$$I_n = [1+ 2+\ldots+n, 1+ 2+\ldots+n+1 )\cap \N \, .
$$

Define two sequences of finite subsets $A_n, B_n$ described by $A_n = \{ r_{2i-1} \mid  i\in I_n\}$ and $B_n =\{ r_{2i} \mid i\in I_n \}\, $.

To simplify notation, in what follows we denote the relator $r_{2i-1}$ by $\alpha_i$ and $r_{2i}$ by $\beta_i$, respectively. Thus, $A_n =\{\alpha_i \mid i\in I_n \}$ and $B_n =\{ \beta_i \mid i\in I_n \}$.

Given a finite subset $X$ in the collection of sets $A_n, B_n$, we construct an element $x\in G$ corresponding to it, as follows. Assume $X$ is composed of the relators $\rho_1,\ldots,\rho_k$ enumerated in increasing order. For every $i\in \{1,2,\ldots,k\}$ let $y_i$ be the prefix of $\rho_i$ of length $\left\lfloor \frac{|\rho_i|}{2} \right\rfloor$.
 Define the element $x= y_1y_2\cdots y_k$. An argument very similar to the one in Lemma \ref{lema:T1ab} implies that $x$ is an $\eta$--relator-tied element and that every geodesic joining $1$ and $x$ is contained in the $\eta$--succession of contours $t_1, y_1t_2,y_1y_2t_3,\ldots,[y_1\cdots y_{k-1}]t_k\, ,$ where $t_i$ is the loop through $1$ in the Cayley graph, labeled by $\rho_i$.

 When $X=A_n$, respectively $X=B_n$ the corresponding element $x$ is denoted by $a_n$, respectively $b_n$.

 We define $\fu_n = [a_n, b_n]$. This implies property \eqref{h2} in Proposition \ref{prop:hn}.

 Lemma \ref{lema:T1ab} applied to geodesics joining $1$ to $\fu_n$ implies that the length $|\fu_n|$ is at least the double of $\left(\frac{1}{2}-\lambda \right)\sum_{i\in I_n} \left[ |\alpha_i| + |\beta_i| \right]\, .$ It follows that property \eqref{h1} in Proposition \ref{prop:hn} is also satisfied.

 We now define the sequence of quasi-morphisms. We start with a general construction. Let $\fv$ be an $\eta$--relator-tied element in $G$.

\begin{defn}
\begin{enumerate}
  \item  Let $(a,b)\in G\times G$. A \emph{quasi-copy of $\fv$ nested inside $(a,b)$} is a pair of points $x,y\in \G^\eta (a,b)$ such that $y=x\fv$ and such that there exists an $\eta$--succession of contours $t_1,\ldots,t_k$ contained in $\G^\eta (a,b)$ such that:
  \begin{itemize}
    \item  $x$ is either one of the endpoints of the intersection of $t_1$ with a contour $t_0$ such that $t_0,t_1,\ldots,t_k$ is an $\eta$--succession contained in $\G^\eta (a,b)$, or the intersection of $t_1$ with a compulsory geodesic preceding  $t_1,\ldots,t_k$ in $\G^\eta (a,b)$;
    \item $y$ is either one of the endpoints of the intersection of $t_k$ with a contour $t_{k+1}$ such that $t_1,\ldots,t_k, t_{k+1}$ is an $\eta$--succession contained in $\G^\eta (a,b)$, or the intersection of $t_k$ with a compulsory geodesic succeeding  $t_1,\ldots,t_k$ in $\G^\eta (a,b)$.
  \end{itemize}

  \me

  \item   We say that two quasi-copies of $\fv$ nested inside $a,b$ are \emph{non-overlapping} if the corresponding $\eta$--successions of contours $t_1,\ldots,t_k$ respectively $\tau_1,\ldots,\tau_k$ are disjoint, as finite sets of contours.

 \item When $(a,b)= (1,g)$ for some element $g\in G$ we speak about \emph{quasi-copies of $\fv$ nested inside $g$}.

\end{enumerate}
\end{defn}

Note that according to the definition of $\G^\eta (a,b)$ and to Lemma \ref{lem:rt2}, the pair of points $x,y$ uniquely determines the $\eta$--succession $t_1,\ldots,t_k$.

 \begin{lem}
Let $x,y$ be a pair of points in $\G^\eta (g)$ composing a nested quasi-copy of $\fv$ in $g$, and let $t_1,\ldots,t_k$ be the corresponding $\eta$--succession of contours. There exists no other pair of points $p,q$ in $\bigcup_{i=1}^k t_i$ such that $q=p\fv\, $.
 \end{lem}

 \proof Lemma \ref{lema:T1ab} can be easily generalized to pairs of points $a,b$ contained in an $\eta$--succession of contours. Applied to the pair $x,y$, it implies that every geodesic joining $x,y$ is $\eta$--relator-tied. This implies that $\fv$ is an $\eta$--relator-tied element. Let $\fg$ be an $\eta$--relator-tied geodesic joining $1$ and $\fv$. It follows that $x\fg$ is contained in $\bigcup_{i=1}^k t_i$, whence the unique sequence of vertices on $\fg$ described in Lemma \ref{lem:rt1} contains $k$ pairs $x_i,y_i$.

 Assume that there exists another pair of points $p\in t_r$ and $q\in t_s$ with $1\lqq r\lqq s\lqq k$ such that $p,q$ compose a nested quasi-copy of $\fv$ in $g$. The $p\fg$ is a geodesic joining $p$ and $q$, which according to Lemma \ref{lema:T1ab} is contained in $\bigcup_{i=r}^s t_i$. The uniqueness of the sequence in Lemma \ref{lem:rt1} implies that $s-r+1=k$, whence $r=1$ and $s=k$. The same uniqueness implies that each pair $px_i, py_i$,  translate of the corresponding pair on $\fg$, is the pair of endpoints of the intersection $p\fg \cap t_i$.

 The first pair in the unique sequence of vertices on $\fg$ as in Lemma \ref{lem:rt1} is of the form $1, h$, where $h$ is represented by a word $w_1$ in $S^{\frac{1}{2}-2\lambda }(R)$, prefix of a relator $\rho$ labeling a unique loop $\tau$ through $1$ in the Cayley graph. By the above $x\tau = p\tau = t_1$, therefore $p\iv x \tau = \tau\, $. This and the small cancelation condition $C' \left(1/12 \right)$ imply that the element $p\iv x$ is trivial in $G$. Indeed, the condition $C' \left(1/12 \right)$ implies that the stabilizer in $G$ of any contour is trivial, otherwise one could find two distinct copies of the same long sub-word in the label of that contour.

 We conclude that $p=x$, and $q= p\fv = x\fv =y\, $.\endproof

 \begin{defn}
 The point $x$ is called the \emph{initial point} of the nested quasi-copy, while $y$ is called the \emph{terminal point} of the nested quasi-copy.
 \end{defn}

 We define $c_\fv :G \times G\to \R$ such that $c_\fv (a,b)$ is the maximal number of pairwise non-overlapping quasi-copies of $\fv $ nested inside $(a,b)$.

By abuse of notation, we define $c_\fv :G \to \R$ such that $c_\fv (g)$ is the maximal number of pairwise non-overlapping quasi-copies of $\fv $ nested inside $g$.

Clearly $c_\fv (a,b) = c_\fv (ha,hb)$ and $c_\fv (g)=  c_\fv (h,hg)$, for every $h\in G$.

 \begin{prop}\label{prop:h3}
 Let $\fv$ be one of the elements $\fu_n$ for $n\in \N$.
 The map $\fh_\fv :G \to \R$, $\fh_\fv = c_\fv-c_{\fv\iv} $ is a quasi-morphism with defect at most $2$.
 \end{prop}

\proof Let $g$ and $h$ be two arbitrary elements in $G$. Our goal is to show that
$$
\left| \fh_\fv (gh)- \fh_\fv(g)-\fh_\fv (h) \right|\leqslant2\, .
$$

The study of geodesic triangles that was done in the preceding section implies that the intersection $\G^\eta (g) \cap \G^\eta (h) \cap \G^\eta (g, gh)$ is either a contour, or a tripod (with some branches possibly reduced to a point) appearing as intersection of three contours, or a sub-path in a contour $\omega$ composed of three consecutive sub-paths (possibly reduced to a point) of lengths $< \lambda |\omega|$, for the first and third, and $<\eta |\omega |$ for the second.  Note that whatever the geometric nature of the intersection, it splits each of the three sets $\G^\eta (g),\, \G^\eta (h),\,  \G^\eta (g, gh),\, $ into two connected components.

We call the intersection $\G^\eta (g) \cap \G^\eta (h) \cap \G^\eta (g, gh)$ the \emph{median object} for the triple $g,h,gh$, and we denote it $\fm (g,h)$.

We say that $\fm (g,h)$ \emph{separates a quasi-copy of $\fv$ nested inside $(a,b)$}, where $(a,b)\in \{(1,g), (1,gh), (g,gh) \}$ if the two points $x,y$ determining that quasi-copy are in two different connected components of $\G^\eta (a,b)\setminus \fm (g,h)$.

Assume that the maxima $c_{\fv^{\pm 1}} (g), \, c_{\fv^{\pm 1}} (gh)$ and $c_{\fv^{\pm 1}} (g, gh)$ are all attained only by considering nested quasi-copies that are not separated by $\fm (g,h)$. In that case one can easily see that
$\fh_\fv (gh)- \fh_\fv(g)-\fh_\fv (h)=0$.

Assume now that every counting that realizes the maximum $c_{\fv } (gh)$ must take into account a pair $x,y$ separated by $\fm (g,h)\, $. Inside $\G^\eta (gh)$ one has then an $\eta$--succession of contours $t_1,\ldots,t_k$ with $x\in t_1$ and $y\in t_k$. The choice of the labels of contours in $\G^\eta (\fu_n )$ implies that:

 \begin{itemize}
   \item no quasi-copy of $\fv\iv$ nested inside $gh$ can contain a sub-sequence in the sequence of contours $t_1,\ldots,t_k$;
   \item  no initial point of a quasi-copy of $\fv$ nested inside $g$ can be contained in $\bigcup_{i=1}^k t_i\cap \G^\eta (g)$;
   \item no terminal point of a quasi-copy of $\fv$ nested inside $(g,gh)$ can be contained in $\bigcup_{i=1}^k t_i \cap \G^\eta (g,gh)$.
 \end{itemize}

It is nevertheless possible that $\bigcup_{i=1}^k t_i\cap \G^\eta (g)$ contains an initial point of a quasi-copy of $\fv\iv$ nested inside $g$. But in that case no terminal point of a quasi-copy of $\fv\iv$ nested inside $(g,gh)$ can be contained in $\bigcup_{i=1}^k t_i \cap \G^\eta (g,gh)$. We thus obtain that
\begin{equation}\label{eq:g2}
\fh_\fv (gh)- \fh_\fv(g)-\fh_\fv (h)=2\, .
\end{equation}

Similarly, $\bigcup_{i=1}^k t_i \cap \G^\eta (g,gh)$ may contain a terminal point of a quasi-copy of $\fv\iv$ nested inside $(g,gh)$; in which case $\bigcup_{i=1}^k t_i\cap \G^\eta (g)$ cannot contain an initial point of a quasi-copy of $\fv\iv$ nested inside $g$, and \eqref{eq:g2} is still verified.

If none of the above two cases occurs then the right-hand side in \eqref{eq:g2} is $1$.

In the case when every counting that realizes the maximum $c_{\fv\iv } (gh)$ must take into account a pair $x,y$ separated by $\fm (g,h)\, $ similar arguments work and give equalities like in  \eqref{eq:g2}, with the right hand side either $-2$ or $-1$.

The cases when $c_{\fv^{\pm 1} } (gh)$ is replaced by either $c_{\fv^{\pm 1}} (g)$ or  $c_{\fv^{\pm 1} } (g,gh)$ are treated similarly and give equalities like in  \eqref{eq:g2}, with the right hand side $\pm 2$ or $\pm 1$.\endproof

\me

We now finish the proof of Proposition \ref{prop:hn}. Proposition \ref{prop:h3} implies that all the quasi-homomorphisms $\fh_{\fu_n}$ has defect bounded by $2$. Properties \eqref{h4} and \eqref{h5} follow from Corollary \ref{cor:seqr} and from the construction of the $\eta$--relator-tied elements $\fu_n$. \endproof

\me

The end of the proof now follows the standard argument in the work of Epstein-Fujiwara \cite{Epstein-Fujiwara, Fujiwara:BC1, Fujiwara:BC2} and Bestvina-Fujiwara \cite{BestvinaFujiwara}. We repeat it here for the sake of completeness.

\me

\begin{thm}\label{thm:qm}
Let $G$ be an infinitely presented finitely generated group given by a presentation satisfying the small cancelation condition $C'(1/12)$.
Then there exists an injective linear map $\ell^1 \to \widetilde{QH} (G)\, $.
In particular, the dimension of $\widetilde{QH} (G)$ is power continuum.
\end{thm}

\proof We consider the map $\ell^1 \to \mathcal{Q}(G)$ defined by $(a_n) \mapsto \sum_{n} a_n \fh_{\fu_n}\, $. Proposition~\ref{prop:hn}, \eqref{h3}, implies that each image is indeed a quasi-morphism. Proposition \ref{prop:hn}, \eqref{h1}, implies that when $a_n \fh_{\fu_n}$ is evaluated in some element $g\in G$, only finitely many terms take non-zero value, thus the sum is always finite.

The above map defines a linear map $\ell^1 \to \widetilde{QH} (G)\, $. We now prove that it is injective. Let $(a_n)\in \ell^1$ be such that $\fh = \sum_{n} a_n \fh_{\fu_n}\, $ is at bounded distance from a homomorphism. In particular, it follows by Proposition \ref{prop:hn}, \eqref{h2}, that for every $n$ and $k$, $\fh \left(\fu_n^k \right)$ is uniformly bounded.

On the other hand, Proposition \ref{prop:hn}, \eqref{h4} and \eqref{h5}, imply that $\fh\left(\fu_n^k \right)=k$. This gives a contradiction.\endproof

\newcommand{\etalchar}[1]{$^{#1}$}
\def\cprime{$'$} \def\cprime{$'$} \def\cprime{$'$}
\providecommand{\bysame}{\leavevmode\hbox to3em{\hrulefill}\thinspace}
\providecommand{\MR}{\relax\ifhmode\unskip\space\fi MR }
\providecommand{\MRhref}[2]{%
  \href{http://www.ams.org/mathscinet-getitem?mr=#1}{#2}
}
\providecommand{\href}[2]{#2}

\end{document}